\documentclass[a4paper,12pt]{article}
\usepackage{authblk}
\usepackage{amsmath}
\usepackage{amsfonts}
\usepackage{amssymb}
\usepackage{amsthm}
\usepackage{graphicx}
\newtheorem{theorem}{Theorem}

\newtheorem{remark}[theorem]{Remark}

\theoremstyle{remark}

\newcommand{\id}{ {1\!\!\!\:1 } }

\begin{document}

\title{Existence of minimizers in the geometrically non-linear 6-parameter resultant shell
theory with drilling rotations}

\author{Mircea B\^{\i}rsan%
\thanks{\,Mircea B\^{\i}rsan,  Fakult\"at f\"ur  Mathematik, Universit\"at Duisburg-Essen, Campus Essen, Thea-Leymann Str. 9, 45127 Essen, Germany, email: mircea.birsan@uni-due.de ; and
Department of Mathematics, University ``A.I. Cuza'' of Ia\c{s}i, 700506 Ia\c{s}i,  Romania}\,\,\,
and\, Patrizio Neff\,%
\thanks{Corresponding author: Patrizio Neff, Lehrstuhl f\"ur Nichtlineare Analysis und Modellierung, Fakult\"at f\"ur  Mathematik, Universit\"at Duisburg-Essen, Campus Essen, Thea-Leymann Str. 9, 45127 Essen, Germany, email: patrizio.neff@uni-due.de, Tel.: +49-201-183-4243}
}

% \date{Received: date / Accepted: date}
% The correct dates will be entered by the editor

\maketitle

\thanks{\small{Dedicated to W. Pietraszkiewicz, our friend and source of inspiration for the 6-parameter resultant shell model.}}

\begin{abstract}
The paper is concerned with the geometrically non-linear theory of 6--parame-tric elastic shells with drilling degrees of freedom. This theory establishes a general model for shells, which is characterized by two independent kinematic fields: the translation vector and the  rotation tensor. Thus, the kinematical structure of 6-parameter shells is identical to that of Cosserat shells. We show the existence of global minimizers for the geometrically non-linear 2D equations of elastic shells. The proof of the existence theorem is based on the direct methods of the calculus of variations using essentially the convexity of the energy in the strain and curvature measures. Since our result is valid for general anisotropic shells, we analyze separately the particular cases of isotropic shells, orthotropic shells, and composite shells.
\end{abstract}

\smallskip

\noindent \textbf{Mathematics Subject Classification  (MSC 2010):}  74K25, 74G65, 74G25, 74E10.

\section{Introduction}\label{sect1}

In recent years there has been a revived interest in 2D shell models because of unconventional materials and extremely small aspect to thickness ratio, such as for instance thin polymeric films or biological membranes. For classical engineering materials and for non-extreme aspect to thickness ratios, available 3D FEM-Codes may readily be used such that the need for a truly 2D shell model does not arise anymore. However, for ultrathin specimens the application of a 3D  constitutive law is not  clear at all.  In these extreme cases one is led to employ a 2D shell model.
This paper is concerned with one such model, the geometrically non-linear resultant theory of shells. We consider the 6-parameter model of shells which involves two independent kinematic fields: the translation vector field and the rotation tensor field (in total 6 independent scalar kinematic variables). This theory of shells is one of the most general, and it is also very effective in the treatment of complex shell problems, as can be seen from the works \cite{Pietraszkiewicz10,Eremeyev11,Pietraszkiewicz11}, among others. The resultant  6-parameter theory of shells was originally proposed by Reissner \cite{Reissner74} and it has been considerably developed subsequently. An account of these developments and main achievements have been presented in the books of Libai and Simmonds \cite{Libai98} and    Chr\'o\'scielewski, Makowski and  Pietraszkiewicz \cite{Pietraszkiewicz-book04}.
In this approach,
the 2D equilibrium equations and static boundary conditions of the shell are derived exactly
by direct through-the-thickness integration of the stresses in the 3D balance laws of linear and angular momentum. The kinematic fields are then constructed on the 2D level using the integral identity of the virtual work principle. Following this procedure, the 2D model is expressed in terms of stress resultants and work--averaged deformation fields defined on the shell base surface. It is interesting  that the kinematical structure of 6-parameter shells (involving the translation vector and rotation tensor) is identical to the kinematical structure of Cosserat shells (defined as material surfaces endowed with a triad of rigid directors describing the orientation of points). From this point of view, the 6-parameter theory of shells is related to the shell model proposed initially by the Cosserat brothers \cite{Cosserat09neu} and developed by many authors, such as  Zhilin \cite{Zhilin76}, Zubov \cite{Zubov97}, Altenbach and Zhilin \cite{Altenbach04}, Eremeyev and Zubov \cite{Eremeyev08}, B\^{\i}rsan and Altenbach \cite{Birsan10}. Using the so--called derivation approach, Neff \cite{Neff_plate04_cmt,Neff_plate07_m3as} has established independently a Cosserat--type model for initially planar shells (plates) which is very similar to the  6-parameter resultant shell model. A   comparison between these two models has been presented in the paper \cite{Birsan-Neff-Plates}, in the case of plates.

On the other hand, we should mention that the kinematic structure of the 6-parameter shell model is different from the kinematic structure of the so--called \emph{Cosserat surfaces}, which are defined as material surfaces with one or more deformable directors attached to every point, see \cite{Naghdi72,Antman95,Rubin00,Rubin04,Altenbach-Erem12,Altenbach-Erem-Review}. For instance, the kinematics of Cosserat surfaces with one deformable director  is characterized also by 6 degrees of freedom (3 for the position of material points and 3 for the orientation and stretch of the material line element through the thickness), which differ essentially from the 6 degrees of freedom in the 6-parameter resultant shell model.

The topic of existence of solutions for the 2D equations of linear and non-linear elastic shells has been treated in many works.   The results that can be found in the literature refer to various types of shell models and they employ different techniques, see e.g. \cite{Koiter60,Simo89.1,Simo92,Steigmann08,Steigmann12,Sansour92,Aganovic06,Aganovic07,Tiba02,Davini75,Birsan08,Wisniewski10,Wisniewski12,Ibrahim94,Badur-Pietrasz86}. The method of formal asymptotic expansions is one method of investigation which allows   for the derivation and justification of plate and shell models. The existence theory for linear or nonlinear shells is presented in details in the books of Ciarlet \cite{Ciarlet97,Ciarlet98,Ciarlet00}, together with many historical remarks and bibliographic references. Another fruitful approach to the existence theory of 2D plate and shell models (obtained as limit cases of 3D models) is the $\Gamma$-convergence analysis of thin structures, see e.g. \cite{Neff_Chelminski_ifb07,Neff_Danzig05,Neff_Hong_Reissner08,Paroni06,Paroni06b}. Concerning the geometrically  non-linear 6-parameter theory of elastic shells, there is no existence theorem published in the literature yet, as far as we are aware of. In the case of \emph{linear} micropolar shells, the existence of weak solutions has been recently proved in \cite{EremeyevLebedev11}. Existence results for the related (very similar) Cosserat--type model of initially planar shells have been established by Neff \cite{Neff_plate04_cmt,Neff_plate07_m3as}. In \cite{Neff_Chelminski_ifb07,Neff_Danzig05,Neff_Hong_Reissner08} the linearized version of this model has been analyzed and compared with the
classical membrane and bending plate models given by the Reissner--Mindlin or Kirchhoff--Love theories.

In the present work, we prove the existence of minimizers for the minimization problem of the total potential energy associated to the deformation of geometrically non-linear 6-parameter elastic shells. We emphasize that our work is not concerned with the derivation of the 2D shell model, but it presents existence results for the well-established 2D theory of 6-parameter elastic shells. It should be mentioned from the beginning that this model refers to shells made of a simple (classical) elastic material, not a generalized (Cosserat or micropolar) continuum. However, the rotation tensor field appears naturally in this theory, in the course of the exact through-the-thickness reduction of the 3D formulation of the problem to the 2D one \cite{Libai98,Pietraszkiewicz-book04,Eremeyev06}. Thus, in spite of the above mentioned similarity with the kinematics of Cosserat shells, the material of the shell in the resultant 6-parameter model is described as a simple continuum (without any specific microstructure or material length scale). On the other hand, in the case of dimensional reduction of the 3D equations of micropolar shell-like bodies one can obtain the same 6-parameter theory with modified 2D constitutive equations, see e.g. \cite{Altenbach-Erem09} for the linear case and \cite{Neff_plate04_cmt,Neff_plate07_m3as,Zubov09} for the nonlinear case, or one can obtain more complex theories as in \cite{Eringen67}.

For the proof of existence, we employ the direct methods of the calculus of variations and extend the techniques presented in  \cite{Neff_plate04_cmt,Neff_plate07_m3as} to the case of general shells (with non-vanishing curvature in the reference configuration). In Section 2 we present briefly the kinematics of general  6-parameter shells and the equations of equilibrium. In Section 3 we give some alternative formulas for the strain tensor and  curvature tensor, which are written in the direct tensorial notation as well as in the component (matrix) notation. These expressions are needed subsequently in the proof of our main result.
In Section 4 we formulate the two-field minimization problem for general elastic shells, corresponding to mixed--type boundary conditions. Under the assumptions of convexity and coercivity of the quadratic strain energy function (physically linear material response), we prove the existence of minimizers over a large  set of admissible pairs. Thus, the minimizing solution pair is of class $\boldsymbol{H}^1(\omega, \mathbb{R}^3)$ for the translation vector and   $\boldsymbol{H}^1(\omega, SO(3))$ for the rotation tensor. The existence result is valid for general anisotropic elastic shells having arbitrary geometry of the reference configuration.
Section 5 includes some applications of the existence theorem and discussions of special cases. We present a convenient way to choose the initial directors and the parametrization of the reference surface. Then, we consider separately the cases of isotropic shells, orthotropic shells, and composite layered shells and we present   the respective forms of the strain energy densities. Applying the theorem stated previously, we establish the conditions on the constitutive coefficients that ensure the existence of minimizers in each situation. This analysis shows the usefulness of our theoretical result in the treatment of practical  problems for elastic shells.

\section{General 6-parameter resultant shells}\label{sect2}

\begin{figure}
\begin{center}
\includegraphics[scale=1]{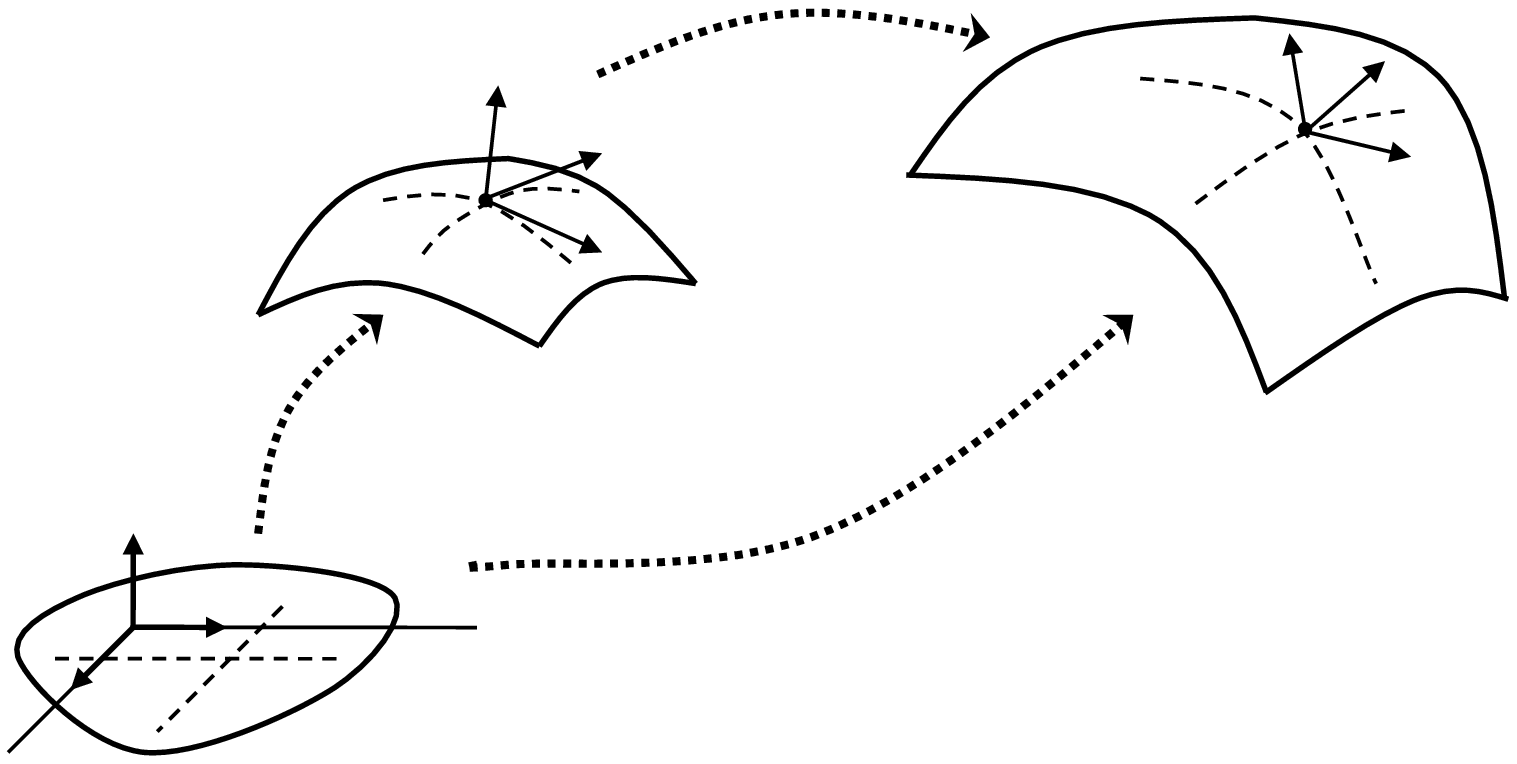}
\put(-398,43){\small{$O$}} \put(-411,20){\small{$\boldsymbol{e}_1$}} \put(-378,45){\small{$\boldsymbol{e}_2$}}
\put(-432,2){$x_1$} \put(-316,31){$x_2$}
    \put(-415,64){\small{$\boldsymbol{e}_3$}} \put(-372,14){$\omega$}
    \put(-357,88){$\boldsymbol{y}^0(x_1,x_2)$} \put(-228,82){$\boldsymbol{R}(x_1,x_2)$}
    \put(-195,60){$\boldsymbol{y}(x_1,x_2)$} \put(-410,105){$\boldsymbol{Q}^{0}(x_1,x_2)$}
    \put(-293,130){${S}^0$} \put(-264.5,145){\small{$\boldsymbol{d}^0_1$}} \put(-264,173){\small{$\boldsymbol{d}^0_2$}} \put(-311,191){\small{$\boldsymbol{d}^0_3$}}
    \put(-67.5,173.5){\small{$\boldsymbol{y} $}} \put(-33,167){$\boldsymbol{d}_1$} \put(-39.5,195){\small{$\boldsymbol{d}_2$}} \put(-81,202){\small{$\boldsymbol{d}_3$}}
\put(-220,204){$\boldsymbol{\chi}( \boldsymbol{y}^0)$}
\put(-220,225){$\boldsymbol{Q }^e( \boldsymbol{y}^0)$} \put(-307,150){\small{$\boldsymbol{y}^0$}}
 \put(-72,120){$ {S}$}
\caption{The base surface $S^0$ of the shell in the initial configuration, the base surface $S$ in the deformed configuration, and the fictitious planar reference configuration $\omega$. The orthonormal triads of vectors $\{\boldsymbol{e}_i\}\,$, $\{\boldsymbol{d}_i^0 \}$ and $\{\boldsymbol{d}_i\}$ are related through the relations $ \boldsymbol{d}_i=\boldsymbol{Q }^e\boldsymbol{d}_i^0=\boldsymbol{R}\boldsymbol{e}_i \,$ and $\boldsymbol{d}_i^0=\boldsymbol{Q}^{0}\boldsymbol{e}_i\,$, where $\boldsymbol{Q}^e  $ is the elastic rotation field, $\boldsymbol{Q}^{0}$ is the  initial rotation, and $\boldsymbol{R}$ is the total rotation field.}
\label{Fig1}       % Give a unique label
\end{center}
\end{figure}

Consider a general 6-parameter shell and denote with $S^0$ the base surface of the shell in the reference (initial) configuration and with $S$ the base surface in the deformed configuration. Let $O$ be a fixed point in the Euclidean space and $\{\boldsymbol{e}_1,\boldsymbol{e}_2,\boldsymbol{e}_3\}$ the fixed orthonormal vector basis. The reference configuration is represented by the position vector $\boldsymbol{y}^0$ (relative to the point $O$) of the base surface $S^0$ \emph{plus} the structure tensor $\boldsymbol{Q}^{0}$. The structure tensor is a second order proper orthogonal tensor which can be described by an orthonormal triad of directors  $\{\boldsymbol{d}^0_1,\boldsymbol{d}^0_2,\boldsymbol{d}^0_3\}$  attached to every point \cite{Pietraszkiewicz-book04,Eremeyev06}. Thus the reference (initial) configuration is characterized by the functions
\begin{equation}\label{1}
\begin{array}{l}
    \,\boldsymbol{y}^0:\omega\subset \mathbb{R}^2\rightarrow\mathbb{R}^3,\qquad\qquad\,\,\boldsymbol{y}^0=\boldsymbol{y}^0(x_1,x_2) ,\\
    \boldsymbol{Q}^{0}:\omega\subset \mathbb{R}^2\rightarrow SO(3), \qquad\,\,\, \boldsymbol{Q}^{0} = \boldsymbol{d}_i^0 (x_1,x_2)\otimes \boldsymbol{e}_i\,,
    \end{array}
\end{equation}
where thus $(x_1,x_2)$ are material curvilinear coordinates on the surface $S^0\,$. Throughout the paper Latin indexes $i,j,...$ take the values $\{1,2,3\}$, while Greek indexes $\alpha,\beta,...$ the values $\{1,2\}$. The usual Einstein summation convention over repeated indexes is employed. We assume that the curvilinear coordinates $(x_1,x_2)\in\omega$  range over a bounded open domain $\omega$ (with Lipschitz boundary $\partial\omega$) of the $Ox_1x_2$ plane, see Figure \ref{Fig1}.
Let us denote the partial derivative with respect to $x_\alpha$ by $\partial_\alpha f=\frac{\partial f}{\partial x_\alpha}\,\,$, for any function $f$. We designate by $\{\boldsymbol{a}_1,\boldsymbol{a}_2 \}$ the (covariant) base vectors in the tangent plane of $S^0$ and by $\boldsymbol{n}^0$ the unit normal to $S^0$ given by
\begin{equation}\label{2}
   \boldsymbol{a}_\alpha= \partial_\alpha \boldsymbol{y}^0 = \dfrac{ \partial\boldsymbol{y}^0}{\partial x_\alpha}\,,\qquad \boldsymbol{n}^0=  \dfrac{ \boldsymbol{a}_1\times\boldsymbol{a}_2}{ \| \boldsymbol{a}_1\times\boldsymbol{a}_2\|}\,\,.
\end{equation}
The reciprocal (contravariant) basis $\{\boldsymbol{a}^1,\boldsymbol{a}^2 \}$ of the tangent plane is defined by
$ \boldsymbol{a}^\alpha\cdot\boldsymbol{a}_\beta =\delta^\alpha_\beta$ (the Kronecker symbol). We also use the notations
$$\boldsymbol{a}_3=\boldsymbol{a}^3=\boldsymbol{n}^0, \quad a_{\alpha\beta}=\boldsymbol{a}_\alpha\cdot \boldsymbol{a}_\beta\,,\quad a^{\alpha\beta}=\boldsymbol{a}^\alpha\cdot\boldsymbol{a}^\beta,\quad a=\sqrt{\det(a_{\alpha\beta})_{2\times 2}}>0.$$

\begin{figure}
\begin{center}
\includegraphics[scale=0.9]{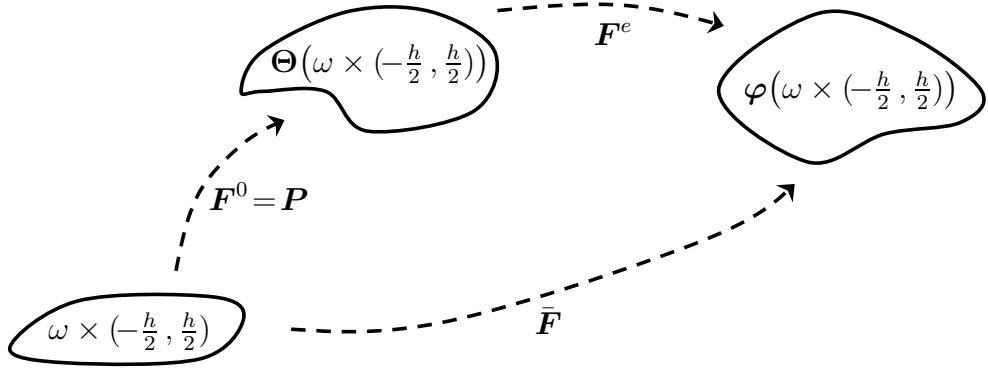}
     \put(-354,14){$\omega\times (\!-\frac{h}{2}\,,\frac{h}{2}) $}
    \put(-271,115){$ \boldsymbol{\Theta}\big( \omega\times(\!-\frac{h}{2}\,,\frac{h}{2})\big) $}
 \put(-94,105){$ \boldsymbol{\varphi}\big( \omega\times(\!-\frac{h}{2}\,,\frac{h}{2})\big) $}
 \put(-294,63){$\boldsymbol{F}^0\!=\!\boldsymbol{P} $}
    \put(-172,16){$ \bar{\boldsymbol{F}}$}
 \put(-150,125){$ \boldsymbol{F}^e$}
\caption{\small{Multiplicative decomposition of the total deformation gradient $\,\bar{F}\,=\,F^eF^0\,$ into the elastic shell deformation gradient $F^e$ and the initial deformation gradient $F^0=P$. Interpretation in terms of reconstructed three-dimensional quantities. The elastic response is governed by $F^e$. The curved initial configuration corresponds to the intermediate stress free configuration in multiplicative plasticity.}}
\label{Fig2}       % Give a unique label
\end{center}
\end{figure}

For the deformed configuration of the shell, let $\boldsymbol{y}(x_1,x_2)$ denote the position vector (relative to $O$) and $\{ \boldsymbol{d}_i(x_1,x_2)\}$ the orthonormal triad of directors attached to the point with initial curvilinear coordinates $(x_1,x_2)$. The deformed configuration is completely characterized by the functions
\begin{equation}\label{3}
    \boldsymbol{y}=\boldsymbol{\chi}(\boldsymbol{y}^0),\qquad \boldsymbol{Q}^e= \boldsymbol{d}_i\otimes \boldsymbol{d}_i^0\in SO(3),
\end{equation}
where $\boldsymbol{\chi}:S^0\rightarrow\mathbb{R}^3$ represents the deformation of the base surface and the proper orthogonal tensor field $\boldsymbol{Q}^e$ is the (effective)  elastic rotation. The displacement vector is defined as usual by $\boldsymbol{u}= \boldsymbol{y}-\boldsymbol{y}^0\,$.

We mention that the role of the triads of directors $\{\boldsymbol{d}_i^0 \}$ and $\{\boldsymbol{d}_i\}$ is to determine the structure tensor $\boldsymbol{Q}^{0}$ and the rotation tensor $\boldsymbol{Q}^{e}$ of the shell, respectively. Thus, the directors do not describe here any microstructure of the material. According to the derivation procedure of the 6-parameter shell model, the kinematical fields $\boldsymbol{y}$ and $\boldsymbol{Q}^e$ are uniquely defined as the work--conjugate averages of 3D deformation distribution over the shell thickness, whose virtual values enter the virtual work principle of the shell (see \cite{Pietraszkiewicz-book04,Pietraszkiewicz-book02}).

In view of \eqref{1} and \eqref{3}, the deformed configuration can alternatively be characterized by the functions
$$\boldsymbol{y}=\boldsymbol{y}(x_1,x_2)= \boldsymbol{\chi}\big( \boldsymbol{y}^0(x_1,x_2)\big),\qquad
\boldsymbol{R}(x_1,x_2)=\boldsymbol{Q}^e\boldsymbol{Q}^{0}= \boldsymbol{d}_i(x_1,x_2)\otimes \boldsymbol{e}_i\in SO(3),$$
where the vector $\boldsymbol{y}$ and the orthogonal tensor $\boldsymbol{R}$ are fields defined over $\omega$.  The  orthogonal tensor field   $\boldsymbol{Q}^e$ represents the \emph{elastic rotation} tensor between the reference and deformed configurations \cite{Pimenta04,Pimenta09}. The tensor $\boldsymbol{Q}^{0}$ is the     \emph{initial rotation}  field, while $\boldsymbol{R}= \boldsymbol{Q}^e\boldsymbol{Q}^{0} $ describes the
\emph{total rotation} from the fictitious planar reference configuration $\omega$ (endowed with the triad $\{\boldsymbol{e}_i\}$) to the deformed configuration $S$. We mention that the tensors $\boldsymbol{Q}^{0}$   and $\boldsymbol{R}$ are also called the \emph{structure tensors} of the reference and deformed configurations, respectively \cite{Pietraszkiewicz-book04,Eremeyev06}.
The following relations hold
\begin{equation}\label{4}
    \boldsymbol{Q}^e=\boldsymbol{R}\boldsymbol{Q}^{0,T}\,,\qquad \boldsymbol{d}_i^0=\boldsymbol{Q}^{0} \boldsymbol{e}_i\,,\qquad \boldsymbol{d}_i=\boldsymbol{Q}^e\boldsymbol{d}_i^0= \boldsymbol{R}  \boldsymbol{e}_i\,.
\end{equation}
Usually, the initial directors $\,\boldsymbol{d}_i^0\,$ are chosen such that $\,\boldsymbol{d}_3^0 = \boldsymbol{n}^0$ and $ \boldsymbol{d}_\alpha^0$ belong to the tangent plane of $S^0$ (see Remark \ref{rem6}). This assumption is not necessary in general and we do not use it in the proof of our existence result.

Let $\boldsymbol{F}=\text{Grad}_s\boldsymbol{y}=\partial_\alpha \boldsymbol{y}\otimes\boldsymbol{a}^\alpha$ denote the (total) shell deformation gradient tensor. The strong form of the equations of equilibrium for 6-parameter shells can be written in the form \cite{Eremeyev06}
\begin{equation}\label{5}
    \mathrm{Div}_s\, \boldsymbol{N}+\boldsymbol{f}=\boldsymbol{0},\qquad \mathrm{Div}_s\, \boldsymbol{M} + \mathrm{axl}(\boldsymbol{N}\boldsymbol{F}^T-\boldsymbol{F}\boldsymbol{N}^T)
    +\boldsymbol{c}=\boldsymbol{0},
\end{equation}
where $\boldsymbol{N}$ and $\boldsymbol{M}$ are the internal surface stress resultant and stress couple tensors of the 1$^{st}$ Piola--Kirchhoff type, while $\boldsymbol{f}$ and $\boldsymbol{c}$ are the external surface resultant force and couple vectors applied to points of $S$, but measured per unit area of $S^0\,$. The operators Grad$_s$ and  Div$_s$ are the surface gradient and surface divergence, respectively, defined intrinsically in \cite{Gurtin-Murdoch-75,Murdoch-90}. The superscript $(\cdot)^T$ denotes the transpose and axl$(\,\cdot)$ represents the axial vector of a skew--symmetric tensor.

Let $\boldsymbol{\nu}$ be the external unit normal vector to the boundary curve $\partial S^0$ lying in the tangent plane. We consider boundary conditions of the type \cite{Pietraszkiewicz04,Pietraszkiewicz11}
\begin{equation}\label{6}
\begin{array}{l}
\boldsymbol{N}\boldsymbol{\nu}=\boldsymbol{n}^*,\qquad \boldsymbol{M}\boldsymbol{\nu}=\boldsymbol{m}^*\qquad\mathrm{along}\,\,\,\partial S^0_f\,,\\
    \quad\boldsymbol{y}=\boldsymbol{y}^* ,\qquad\quad\,\, \boldsymbol{R}=\boldsymbol{R}^* \qquad\mathrm{along}\,\,\,\partial S^0_d\,,
    \end{array}
\end{equation}
where $\partial S^0=\partial S^0_f\cup \partial S^0_d $ is a disjoint partition  of $S^0\,$ ($\partial S^0_f\cap \partial S^0_d=\emptyset$) with length($\partial S^0_d )>0$. Here, $\boldsymbol{n}^*$ and $\boldsymbol{m}^*$ are the external boundary resultant force and couple vectors respectively, applied along the deformed boundary $\partial S$, but measured per unit length of $\partial S^0_f\subset\partial S^0\,$. We denote by $\partial \omega_f$ and $\partial \omega_d$ the subsets of the boundary curve $\partial \omega$ which correspond to $\partial S^0_f$ and $\partial S^0_d$ respectively, through the mapping $\boldsymbol{y}^0\,$.

The weak form associated to these local balance equations for shells has been presented in \cite{Libai98,Pietraszkiewicz-book04,Pietraszkiewicz04}.

\section{Elastic shell strain and curvature measures}\label{sect3}

According to \cite{Eremeyev06,Pietraszkiewicz-book04}, the elastic shell strain tensor $\boldsymbol{E}^e$ in the material representation is given by
\begin{equation}\label{7}
    \boldsymbol{E}^e=\boldsymbol{Q}^{e,T}\text{Grad}_s\,\boldsymbol{y}- \text{Grad}_s\,\boldsymbol{y}^0=\big(\boldsymbol{Q}^{e,T}\partial_\alpha \boldsymbol{y}- \partial_\alpha \boldsymbol{y}^0 \big) \otimes \boldsymbol{a}^\alpha\,,
\end{equation}
since $\text{Grad}_s\boldsymbol{y}=\partial_\alpha \boldsymbol{y} \otimes\boldsymbol{a}^\alpha\,$.
It is useful to write  the strain tensor $\boldsymbol{E}^e$ in the alternative form
\begin{equation*}
\begin{array}{l}
\boldsymbol{E}^e=\big(\boldsymbol{Q}^{e,T}\partial_\alpha \boldsymbol{y} - \boldsymbol{a}_\alpha\big) \otimes \boldsymbol{a}^\alpha=  \big(\boldsymbol{Q}^{e,T}\partial_\alpha \boldsymbol{y}\otimes \boldsymbol{a}^\alpha+  \boldsymbol{n}^0  \otimes \boldsymbol{a}^3\big)- \big(  \boldsymbol{a}_i\otimes \boldsymbol{a}^i\big) \\
\quad\,\,\, =  \big(\boldsymbol{Q}^{e,T}\partial_\alpha \boldsymbol{y}\otimes \boldsymbol{e}_\alpha+  \boldsymbol{n}^0  \otimes \boldsymbol{e}_3\big)  \big(  \boldsymbol{e}_i\otimes \boldsymbol{a}^i\big) -\id_3
\end{array}
\end{equation*}
or equivalently, since $\big(  \boldsymbol{e}_i\otimes \boldsymbol{a}^i\big)= \big(  \boldsymbol{a}_i\otimes \boldsymbol{e}_i\big)^{-1}$,
\begin{equation}\label{13}
\begin{array}{l}
    \boldsymbol{E}^e  =  \big(\boldsymbol{Q}^{e,T}\partial_\alpha \boldsymbol{y}\otimes \boldsymbol{e}_\alpha+  \boldsymbol{n}^0  \otimes \boldsymbol{e}_3\big)  \big(  \boldsymbol{a}_i\otimes \boldsymbol{e}_i\big)^{-1} -\id_3 \,\,\,=\,\, \, \overline{\boldsymbol{U}}^{\,e} -\id_3\,\,,\vspace{4pt}\\
    \quad\text{with}\qquad \overline{\boldsymbol{U}}^{\,e}= \big(\boldsymbol{Q}^{e,T}\partial_\alpha \boldsymbol{y}\otimes \boldsymbol{e}_\alpha+  \boldsymbol{n}^0  \otimes \boldsymbol{e}_3\big)  \big(  \boldsymbol{a}_i\otimes \boldsymbol{e}_i\big)^{-1},
    \end{array}
\end{equation}
where $\id_3=\boldsymbol{e}_i\otimes\boldsymbol{e}_i$ is the identity tensor and  $\overline{\boldsymbol{U}}^{\,e}$ represents the non-symmetric \emph{elastic shell stretch tensor}, which can be seen as the 2D analog of the 3D non-symmetric Biot--type stretch tensor \cite{Neff_Biot07} or the first Cosserat deformation tensor \cite[page 123, eq. (43)]{Cosserat09neu} for the shell.
Let us denote by  $\boldsymbol{P}$  the tensor defined by
\begin{equation}\label{14}
    \boldsymbol{P}  =  \boldsymbol{a}_i\otimes \boldsymbol{e}_i = \partial_\alpha \boldsymbol{y}^0\otimes \boldsymbol{e}_\alpha+  \boldsymbol{n}^0  \otimes \boldsymbol{e}_3\,.
\end{equation}
Then, from \eqref{13} and \eqref{14} we get
\begin{equation}\label{14,1}
\begin{array}{l}
    \boldsymbol{E}^e =\, \overline{\boldsymbol{U}}^{\,e} -\id_3\, = \,\boldsymbol{Q}^{e,T} \big( \partial_\alpha \boldsymbol{y}\otimes \boldsymbol{e}_\alpha+  \boldsymbol{Q}^{e}\boldsymbol{n}^0  \otimes \boldsymbol{e}_3\big)  \boldsymbol{P}^{-1}-\id_3\,,\vspace{4pt}\\
    \qquad\,\,\,\,\overline{\boldsymbol{U}}^{\,e}= \,\boldsymbol{Q}^{e,T} \big( \partial_\alpha \boldsymbol{y}\otimes \boldsymbol{e}_\alpha+  \boldsymbol{Q}^{e}\boldsymbol{n}^0  \otimes \boldsymbol{e}_3\big)  \boldsymbol{P}^{-1}.
    \end{array}
\end{equation}
In the sequel, it is useful to write the elastic shell strain tensors   in component form, relative to the fixed tensor basis $\{\boldsymbol{e}_i\otimes\boldsymbol{e}_j\}$. Let $E^e=\big(E^e_{ij}\big)_{3\times 3}$   be the matrix of components for the tensor $\boldsymbol{E}^e=E^e_{ij}\boldsymbol{e}_i\otimes\boldsymbol{e}_j\,$. In general, we decompose any second order tensor $\boldsymbol{T}$ in the form $\boldsymbol{T}=T_{ij}\boldsymbol{e}_i\otimes\boldsymbol{e}_j$ and denote by  $T=\big(T_{ij}\big)_{3\times 3}$ the matrix of components. Also, for any vector
$\boldsymbol{v}=v_{i }\boldsymbol{e}_i$ we designate by $v=\big(v_{i }\big)_{3\times 1}$ the column matrix of components.
\begin{remark}\label{rem1}
The matrix of components $P=\big(P_{ij}\big)_{3\times 3}$ for the tensor defined in \eqref{14} can be specified in terms of its 3 column vectors as follows
\begin{equation}\label{14bis}
    P=\Big(\,\partial_1 y^{0}\,\Big|\, \partial_2 y^{0}\,\Big|\, n^0\,\Big)_{3\times 3}=\Big(\,\nabla y^{0}\,\Big|\, n^0\,\Big)_{3\times 3}= \Big(\,a_{1}\,\Big|\, a_{2}\,\Big|\, n^0\, \Big)_{3\times 3}\,\,\,.
\end{equation}
We mention that the tensor
$\boldsymbol{P}$ introduced in \eqref{14} can be seen as a three-dimensional  (deformation) gradient
\begin{equation}\label{14tris}
\begin{array}{c}
\boldsymbol{P}= \nabla \, \boldsymbol{\Theta}(x_1,x_2,x_3)_{\big|x_3=0}\,\,,\qquad\mathrm{with}\vspace{4pt}\\
 \boldsymbol{\Theta}(x_1,x_2,x_3):=\boldsymbol{y}^0(x_1,x_2)+x_3\,\boldsymbol{n}^0(x_1,x_2),
\end{array}
\end{equation}
and it satisfies  $\det\boldsymbol{P} =  \sqrt{\det a_{\alpha\beta}}  =a>0$, so that the inverse $\boldsymbol{P}^{-1}$ exists.
The mapping $\,\boldsymbol{\Theta}:\,\omega\times\big(\!-\frac{h}{2}\,,\frac{h}{2}\,\big)\rightarrow\mathbb{R}^3$  has been introduced previously in \cite{Ciarlet00,Ciarlet05,Neff_plate04_cmt} and employed for the geometrical description of 3D shells ($h$ denotes the thickness of the shell).
\end{remark}
By virtue of \eqref{13} and \eqref{14}, we obtain the following  matrix form for the strain tensor $\boldsymbol{E}^e$
\begin{equation}\label{15}
    E^e= \Big(\,Q^{e,T}\partial_1 y\,\Big|\, Q^{e,T}\partial_2 y\,\Big|\, n^0\,\Big) \,P^{-1} - \id_3\,,
\end{equation}
where $\id_3=\big(\delta_{ij}\big)_{3\times 3}$ is the unit matrix. The matrix $E^e$ can be written equivalently as
\begin{equation}\label{15,1}
\begin{array}{l}
    E^e \,=\, Q^{e,T}\Big(\,\partial_1 y\,\Big|\,  \partial_2 y\,\Big|\, Q^{e} n^0\,\Big) \,P^{-1} - \id_3\,,\qquad\text{or}\vspace{4pt}\\
    E^e\,=\,\, \overline{U}^{\,e}-\id_3 \,\,= \,\, Q^{e,T}\,F^e - \id_3 \,\, = \,\, Q^{e,T}\bar{F}  \,P^{-1}-  \id_3\,,
    \end{array}
\end{equation}
with
\begin{equation}\label{15,2}
\begin{array}{l}
    \overline{U}^{\,e}= Q^{e,T}\Big(\,\partial_1 y\,\Big|\,  \partial_2 y\,\Big|\, Q^{e} n^0\,\Big) \,P^{-1} = Q^{e,T}\Big(\,\nabla y\,\Big|\, Q^{e} n^0\,\Big) \,P^{-1}\,,\vspace{4pt}\\
    F^e:= \Big(\,\partial_1 y\,\Big|\,  \partial_2 y\,\Big|\, Q^{e} n^0\,\Big) \,P^{-1} = \Big(\,\nabla y\,\Big|\, Q^{e} n^0\,\Big) \big(\nabla \Theta(x_1,x_2,0)\big)^{-1}\,,  \vspace{4pt}\\
    \bar{F}\,\,:= \Big(\,\partial_1 y\,\Big|\,  \partial_2 y\,\Big|\, Q^{e} n^0\,\Big) = \Big(\,\nabla y\,\Big|\, Q^{e} n^0\,\Big) ,\vspace{4pt}\\
    F^0:=P=\Big(\,\partial_1 y^{0}\,\Big|\, \partial_2 y^{0}\,\Big|\, n^0\,\Big)= \Big(\,\nabla y^0\,\Big|\,   n^0\,\Big), \vspace{4pt}\\
    \qquad\qquad\qquad\,\,\,\qquad\bar{F}\,\,\,=\,\,F^e\,\,F^0\,.
    \end{array}
\end{equation}
In order to see a parallel with the classical  multiplicative decomposition  into elastic and plastic parts from finite elasto-plasticity \cite{Neff_Cosserat_plasticity05,HutterSFB02}, we may interpret $F^e$ as an elastic shell mid-surface deformation gradient and $F^0=P$ as an initial deformation gradient. Both are gradients of suitably defined mappings, see Remark \ref{rem2} and  Figure \ref{Fig2}, in contrast to the case of elasto-plasticity. In our context, the elastic material response is defined in terms of the elastic part of the deformation, e.g. $\,E^e\,=\, Q^{e,T}\,F^e - \id_3 \,$, cf. \eqref{15,1}.

\begin{remark}\label{rem2}
Although the resultant shell model is truly a 2D theory, we may always consider artificially reconstructed three-dimensional quantities. In this sense, similar to the  context of Remark \ref{rem1}, the tensor $\,\bar{\boldsymbol{F}}=\partial_\alpha \boldsymbol{y}\otimes \boldsymbol{e}_\alpha+  \boldsymbol{Q}^e\boldsymbol{n}^0  \otimes \boldsymbol{e}_3\,$, which has the components matrix  $\,\bar{F}\,$ is a three-dimensional deformation gradient
\begin{equation}\label{15,3}
\begin{array}{c}
\bar{\boldsymbol{F}}= \nabla \, \boldsymbol{\varphi}(x_1,x_2,x_3)_{\big|x_3=0}\,\,,\qquad\mathrm{with}\vspace{4pt}\\
 \boldsymbol{\varphi}(x_1,x_2,x_3):=\boldsymbol{y}(x_1,x_2)+x_3\,\boldsymbol{Q}^e(x_1,x_2) \boldsymbol{n}^0(x_1,x_2)\\
\qquad\qquad\qquad\qquad\,\,\, =   \boldsymbol{y}(x_1,x_2)+x_3\,\boldsymbol{Q}^e(x_1,x_2)\,      \nabla \, \boldsymbol{\Theta}(x_1,x_2,0)\boldsymbol{e}_3.
\end{array}
\end{equation}
Here, the mapping $\,\boldsymbol{\varphi}:\,\omega\times\big(\!-\frac{h}{2}\,,\frac{h}{2}\,\big)\rightarrow\mathbb{R}^3$ is a 3D deformation of the body, in terms of the given 2D quantities $\boldsymbol{y}(x_1,x_2)$ and $\boldsymbol{Q}^e(x_1,x_2)$. Similarly,
\begin{equation*}
\begin{array}{c}
\boldsymbol{F}^e:= \nabla \, \boldsymbol{\varphi}^e\Big(\boldsymbol{\Theta}(x_1,x_2,x_3)_{\big|x_3=0}\Big)\,\,,\qquad\mathrm{with}\vspace{4pt}\\
 \boldsymbol{\varphi}^e\big(\boldsymbol{\Theta}(x_1,x_2,x_3)\big):=
 \boldsymbol{\varphi}(x_1,x_2,x_3).
\end{array}
\end{equation*}
However, we note that $\bar{\boldsymbol{F}}$ cannot be interpreted as the  3D deformation
gradient of the real 3D shell, because in general the initial normals become
arbitrarily curved after deformation.
\end{remark}

In terms of the total rotation $\boldsymbol{R}$ and the initial rotation $\boldsymbol{Q}^{0}$, the elastic shell strain tensor is expressed by
\begin{equation}\label{8}
    \boldsymbol{E}^e=\boldsymbol{Q}^{0}\big(\boldsymbol{R}^T\partial_\alpha \boldsymbol{y} - \boldsymbol{Q}^{0,T} \partial_\alpha \boldsymbol{y}^0 \big) \otimes \boldsymbol{a}^\alpha\,.
\end{equation}
Then, we have
$$\boldsymbol{E}^e=\boldsymbol{Q}^{0}\big[\big(\boldsymbol{R}^T\partial_\alpha \boldsymbol{y} \!\!-\! \boldsymbol{Q}^{0,T} \partial_\alpha \boldsymbol{y}^0\big)\! \otimes\! \boldsymbol{e}_\alpha\big] \big(  \boldsymbol{e}_i\otimes \boldsymbol{a}^i\big)= \boldsymbol{Q}^{0}\big[\big(\boldsymbol{R}^T\partial_\alpha \boldsymbol{y}\!-\! \boldsymbol{Q}^{0,T} \partial_\alpha \boldsymbol{y}^0\big) \!\otimes\! \boldsymbol{e}_\alpha\big]
 \boldsymbol{P}^{-1}$$
which can be written in matrix form as follows
\begin{equation}\label{16}
    E^e= Q^0\,H \,P^{-1} \quad\text{with}\quad H:=\Big(\,R^T\partial_1 y- Q^{0,T}\partial_1 y^{0}\,\Big|\, R^T\partial_2 y- Q^{0,T}\partial_2 y^{0}\,\Big|\,\, 0\,\,\Big)_{3\times 3}\, .
\end{equation}\smallskip

On the other hand, the elastic shell curvature  tensor $\boldsymbol{K}^e$ in the material description is defined by \cite{Eremeyev06,Pietraszkiewicz-book04}
\begin{equation}\label{9}
     \boldsymbol{K}^e=\big[\boldsymbol{Q}^{e,T}\text{axl}(\partial_\alpha \boldsymbol{R}  \boldsymbol{R}^T)- \text{axl}(\partial_\alpha \boldsymbol{Q}^{0}\boldsymbol{Q}^{0,T})\big] \otimes \boldsymbol{a}^\alpha\,.
\end{equation}
In order to write $\boldsymbol{K}^e$ in a  form more convenient to us, we use relations of the type
\begin{equation}\label{10}
    \tilde{\boldsymbol{Q}}^{ T}\text{axl}(\partial_\alpha\tilde{ \boldsymbol{Q}} \, \tilde{\boldsymbol{Q}}^{ T})=  \text{axl}(\tilde{\boldsymbol{Q}}^{ T} \partial_\alpha \tilde{ \boldsymbol{Q}}  ),\qquad \text{axl}(\tilde{\boldsymbol{Q}} \boldsymbol{A}\tilde{\boldsymbol{Q}}^{ T}) = \tilde{\boldsymbol{Q}}  \,\text{axl}(\boldsymbol{A}),
\end{equation}
which hold true for any   orthogonal tensor $\tilde{\boldsymbol{Q}} \in SO(3)$ and any skew--symmetric tensor $\boldsymbol{A}\in \frak{so}(3)$ (see e.g., \cite{Birsan-Neff-AnnRom12}). Using \eqref{10} in \eqref{9} we can write the elastic  curvature tensor $\boldsymbol{K}^e$ in the equivalent forms
\begin{equation}\label{11}
    \boldsymbol{K}^e= \boldsymbol{Q}^{e,T}\text{axl}(\partial_\alpha \boldsymbol{Q}^e  \boldsymbol{Q}^{e,T}) \otimes \boldsymbol{a}^\alpha= \text{axl}(\boldsymbol{Q}^{e,T} \partial_\alpha \boldsymbol{Q}^e ) \otimes \boldsymbol{a}^\alpha ,
\end{equation}
or,
$$\boldsymbol{K}^e= \big[\, \text{axl}(\boldsymbol{Q}^{e,T} \partial_\alpha \boldsymbol{Q}^e) \otimes \boldsymbol{e}_\alpha\,\big]
 \big(  \boldsymbol{a}_i\otimes \boldsymbol{e}_i\big)^{-1}  = \big[\, \text{axl}(\boldsymbol{Q}^{e,T} \partial_\alpha \boldsymbol{Q}^e) \otimes \boldsymbol{e}_\alpha\,\big]
 \boldsymbol{P}^{-1}.
$$
Then, the matrix of components $K^e=\big(K^e_{ij}\big)_{3\times 3}$ is given by
\begin{equation}\label{17}
    K^e= \Big(\,\,\text{axl}(Q^{e,T}\partial_1 Q^e) \,\,\Big|\,\, \text{axl}(Q^{e,T}\partial_2 Q^e) \,\,\Big|\, \,0\,\,\Big) \,P^{-1}.
\end{equation}

If we express $\boldsymbol{K}^e$ in terms of the total rotation $\boldsymbol{R}$ and the initial rotation $\boldsymbol{Q}^{0}$, we get
\begin{equation}\label{12}
    \boldsymbol{K}^e= \boldsymbol{Q}^{0} \big[ \text{axl}(\boldsymbol{R}^T \partial_\alpha \boldsymbol{R}  )- \text{axl}(\boldsymbol{Q}^{0,T}\partial_\alpha \boldsymbol{Q}^{0})\big] \otimes \boldsymbol{a}^\alpha\,.
\end{equation}
This relation can be written as
\begin{equation}\label{12,1}
\begin{array}{c}
    \boldsymbol{K}^e=\boldsymbol{K} - \boldsymbol{K}^0,\qquad \text{with} \qquad \boldsymbol{K}:= \boldsymbol{Q}^{0}   \text{axl}(\boldsymbol{R}^T \partial_\alpha \boldsymbol{R}  ) \otimes \boldsymbol{a}^\alpha\,,\vspace{4pt} \\
      \boldsymbol{K}^0:= \boldsymbol{Q}^{0}  \text{axl}(\boldsymbol{Q}^{0,T}\partial_\alpha \boldsymbol{Q}^{0})  \otimes \boldsymbol{a}^\alpha\,\,=\,\,\text{axl}(\partial_\alpha \boldsymbol{Q}^{0}\,\boldsymbol{Q}^{0,T}) \otimes \boldsymbol{a}^\alpha\,,
    \end{array}
\end{equation}
where the tensor $\boldsymbol{K}$ is the total  curvature tensor, while $\boldsymbol{K}^0$ is the initial  curvature (or structure curvature tensor of $S^0$). In view of  \eqref{12} and  \eqref{12,1},
the matrix $K^e=\big(K^e_{ij}\big) $ is given by
\begin{equation}\label{18}
\begin{array}{c}
    K^e\,\,=\,\,Q^0\,L\,P^{-1}\,\,=\,\,K-K^0\qquad\text{with}  \vspace{4pt}\\ L:= \Big(\,\,\text{axl}(R^T\partial_1 R)-\text{axl}(Q^{0,T}\partial_1 Q^{0}) \,\,\Big|\,\, \text{axl}(R^T\partial_2 R) -\text{axl}(Q^{0,T}\partial_2 Q^{0})\,\,\Big|\, \,0\,\,\Big)_{3\times 3}\, ,\vspace{4pt}\\
    \,\, K= Q^0\Big(\,\,\text{axl}(R^T\partial_1 R)  \,\,\Big|\,\, \text{axl}(R^T\partial_2 R) \,\,\Big|\, \,0\,\,\Big)P^{-1}, \vspace{4pt}\\
    \qquad K^0= Q^0\Big(\,\, \text{axl}(Q^{0,T}\partial_1 Q^{0}) \,\,\Big|\,\,  \text{axl}(Q^{0,T}\partial_2 Q^{0})\,\,\Big|\, \,0\,\,\Big)P^{-1}.
    \end{array}
\end{equation}
In what follows, we shall use the expressions \eqref{16} and \eqref{18} of the elastic shell strain measures $\boldsymbol{E}^e$ and $\boldsymbol{K}^e$ written with tensor components in the basis $\{\boldsymbol{e}_i\otimes\boldsymbol{e}_j\}$.
\begin{remark}
As expected, the case of zero strain and bending measures corresponds to a rigid body mode of the shell. Indeed, if $\boldsymbol{E}^e=\boldsymbol{0}$ and $\boldsymbol{K}^e=\boldsymbol{0}$, then from \eqref{7} and \eqref{11} we obtain
$$\partial_\alpha\boldsymbol{y}=\boldsymbol{Q}^e\,\partial_\alpha\boldsymbol{y}^0\qquad \mathrm{and}\qquad \partial_\alpha \boldsymbol{Q}^e=\boldsymbol{0}.$$
Hence, it follows that $\boldsymbol{Q}^e$ is constant and
$$\boldsymbol{y}=\boldsymbol{Q}^e \boldsymbol{y}^0+ \boldsymbol{c}\qquad (\boldsymbol{c}=\mathrm{constant}),$$
which means that the shell undergoes a rigid body motion with constant translation $\,\boldsymbol{c}\,$ and constant rotation $\,\boldsymbol{Q}^e$.
\end{remark}
\begin{remark}
In the case when the base surface $S^0$ of the initial configuration of the shell is planar we may assume that $S^0$ coincides with $\,\omega$. In this situation we have $\,\boldsymbol{a}_i=\boldsymbol{e}_i\,$, $\boldsymbol{P}= \id_3\,$, and
the above strain and curvature measures coincide with those defined for the Cosserat model of planar--shells introduced in \cite{Neff_plate04_cmt,Neff_plate07_m3as}.
\end{remark}
\begin{remark}
In view of \eqref{11} or \eqref{17}, the elastic shell  curvature tensor $\,\boldsymbol{K}^e$ is an analog of the second Cosserat deformation tensor in the 3D theory, see the original Cosserats book \cite[page 123, eq. (44)]{Cosserat09neu}.
\end{remark}

\section{Variational formulation for elastic shells}\label{sect4}

Let us denote the strain energy density of the elastic shell by $W=W(\boldsymbol{E}^e,\boldsymbol{K}^e)$. According to the hyperelasticity assumption, the    internal surface stress resultant $\boldsymbol{N}$ and stress couple tensor $\boldsymbol{M}$ are expressed by the constitutive equations in the form
\begin{equation}\label{19}
    \boldsymbol{N}=\boldsymbol{Q}^e\,\dfrac{\partial\, W}{\partial \boldsymbol{E}^e}\,\,,\qquad \boldsymbol{M}=\boldsymbol{Q}^e\,\dfrac{\partial\, W}{\partial \boldsymbol{K}^e}\,\,.
\end{equation}
In this paper we assume that the strain energy density $W$ is a quadratic function of its arguments $\boldsymbol{E}^e$ and $\boldsymbol{K}^e$. Thus, the considered model  is physically linear and geometrically non-linear. The explicit form of the strain energy function $W$ is presented in \cite{Libai98,Eremeyev06} for isotropic, hemitropic or orthotropic elastic shells.
In general, the coefficients of the strain energy function $W$ depend on the structure curvature tensor $\,\boldsymbol{K}^0$, see \cite{Eremeyev06}.
In \cite{Chroscielewski11}, the case of composite (layered) shells is investigated and the expression of the energy density is established. These special cases will be discussed in Section \ref{sect5}.

Consider the usual Lebesgue spaces $\big(L^p(\omega),\|\cdot\|_{L^p(\omega)}\big)$, $p\geq 1$, and Sobolev space $\big(H^1(\omega),\|\cdot\|_{H^1(\omega)}\big)$. We denote by $\boldsymbol{L}^p(\omega,\mathbb{R}^3)$ (respectively $\boldsymbol{H}^1(\omega,\mathbb{R}^3)$) the space of all vector fields $\boldsymbol{v}=v_i\boldsymbol{e}_i$ such that $v_i\in L^p(\omega)$ (respectively $v_i\in H^1(\omega)$). Similarly, we denote the sets
$    \boldsymbol{H}^1(\omega,\mathbb{R}^{3\times 3})=\{\boldsymbol{T}=T_{ij}\boldsymbol{e}_i\otimes\boldsymbol{e}_j \,|\, T_{ij}\in H^1(\omega)\}$ ,   $\boldsymbol{H}^1(\omega,SO(3))= \{\boldsymbol{T}\in \boldsymbol{H}^1(\omega,\mathbb{R}^{3\times 3}) \,|\,\boldsymbol{T}\in SO(3)\}$ , $ \boldsymbol{L}^p(\omega,\mathbb{R}^{3\times 3})=\{\boldsymbol{T}=T_{ij}\boldsymbol{e}_i\otimes\boldsymbol{e}_j \,|\, T_{ij}\in L^p(\omega)\}$ , $ \boldsymbol{L}^p(\omega,SO(3))= \{\boldsymbol{T}\in \boldsymbol{L}^p(\omega,\mathbb{R}^{3\times 3}) \,|\,\boldsymbol{T}\in SO(3)\}$.
The norm of a tensor $\boldsymbol{T}$ is defined by $\|\boldsymbol{T}\|^2=\text{tr}(\boldsymbol{T} \boldsymbol{T}^T)=T_{ij}T_{ij}\,$.

Concerning the boundary-value problem \eqref{5}, \eqref{6}, we assume the existence of a function $\Lambda(\boldsymbol{y},\boldsymbol{R})$ representing the potential of external surface loads $\boldsymbol{f}$, $\boldsymbol{c}$, and boundary loads $\boldsymbol{n}^*$, $\boldsymbol{m}^*$ (cf. \cite{Pietraszkiewicz04}).

We consider the following two--field minimization problem associated to the deformation of elastic shells: find the pair $(\hat{\boldsymbol{y}},\hat{\boldsymbol{R}})$ in the admissible set $\mathcal{A}$ which realizes the minimum of the functional
\begin{equation}\label{20}
I(\boldsymbol{y},\boldsymbol{R})=\int_{S^0} W(\boldsymbol{E}^e,\boldsymbol{K}^e)\,\mathrm{d}S - \Lambda(\boldsymbol{y},\boldsymbol{R})\qquad\mathrm{for}\qquad (\boldsymbol{y},\boldsymbol{R})\in \mathcal{A},
\end{equation}
where d$S$ is the area element of the surface $S^0\,$. The admissible set $\mathcal{A}$ is defined by
\begin{equation}\label{21}
    \mathcal{A}=\big\{(\boldsymbol{y},\boldsymbol{R})\in\boldsymbol{H}^1(\omega, \mathbb{R}^3)\times\boldsymbol{H}^1(\omega, SO(3))\,\,\big|\,\,\,  \boldsymbol{y}_{\big| \partial S^0_d}=\boldsymbol{y}^*, \,\,\boldsymbol{R}_{\big| \partial S^0_d}=\boldsymbol{R}^* \big\},
\end{equation}
where the boundary conditions are to be understood in the sense of traces.
The tensors $\boldsymbol{E}^e$ and $\boldsymbol{K}^e$ are expressed in terms of $(\boldsymbol{y},\boldsymbol{R})$ through the relations \eqref{8} and \eqref{12}. If we write $W=W(\boldsymbol{E}^e,\boldsymbol{K}^e)=\tilde{W}(\nabla\boldsymbol{y},\boldsymbol{R},\nabla\boldsymbol{R})$, then referring the integral to the (fictitious reference) domain $\omega$, the change of variable formula clearly  gives
\begin{equation}\label{22}
\begin{array}{l}
    \displaystyle{\int_{S^0}} W(\boldsymbol{E}^e,\boldsymbol{K}^e)\,\mathrm{d}S=\displaystyle{\int_\omega } W(\boldsymbol{E}^e,\boldsymbol{K}^e)\,a(x_1,x_2)\, \mathrm{d}x_1\mathrm{d}x_2 \\
    \qquad\,\,\, = \displaystyle{\int_\omega } \tilde{W}\big(\nabla\boldsymbol{y}(x_1,x_2),\boldsymbol{R}(x_1,x_2),\nabla\boldsymbol{R}(x_1,x_2)\big)\, \det\big(\nabla \boldsymbol{\Theta} (x_1,x_2,0)\big)\, \mathrm{d}x_1\mathrm{d}x_2 ,
    \end{array}
\end{equation}
where $a=\sqrt{\det(a_{\alpha\beta})}$ is the notation introduced previously. The variational principle for the total energy of elastic shells with respect to the functional \eqref{20} has been presented in \cite{Pietraszkiewicz04}, Sect.2. We decompose the loading potential $\Lambda(\boldsymbol{y},\boldsymbol{R})$ additively as follows
\begin{equation}\label{23}
\begin{array}{c}
    \Lambda(\boldsymbol{y},\boldsymbol{R})=\Lambda_{S^0}(\boldsymbol{y},\boldsymbol{R}) + \Lambda_{\partial S^0_f}(\boldsymbol{y},\boldsymbol{R}),\\
    \Lambda_{S^0}(\boldsymbol{y},\boldsymbol{R})= \displaystyle{\int_{S^0}}\! \boldsymbol{f}\!\cdot\! \boldsymbol{u}\, \mathrm{d}S + \Pi_{S^0}(\boldsymbol{R}), \quad
    \Lambda_{\partial S^0_f}(\boldsymbol{y},\boldsymbol{R})= \displaystyle{\int_{\partial S^0_f}}\! \boldsymbol{n}^*\!\cdot\! \boldsymbol{u}\, \mathrm{d}l + \Pi_{\partial S^0_f}(\boldsymbol{R}).
    \end{array}
\end{equation}
where $\boldsymbol{u}=\boldsymbol{y}-\boldsymbol{y}^0$ is the displacement vector and d$l$ is the element of length along the boundary curve $\partial S^0_f\,$. In \eqref{23}, $\Lambda_{S^0}(\boldsymbol{y},\boldsymbol{R})$ is the potential of the external surface loads $\boldsymbol{f}, \boldsymbol{c}$, while $\Lambda_{\partial S^0_f}(\boldsymbol{y},\boldsymbol{Q}^e)$ is the potential of the external boundary loads $\boldsymbol{n}^*, \boldsymbol{m}^*$. The expression of the load potential functions $\,\,\Pi_{S^0}\,,  \,\Pi_{\partial S^0_f}:\boldsymbol{L}^2( \omega ,SO(3))\rightarrow\mathbb{R}$ are not given explicitly, but they are assumed to be continuous and bounded operators. Of course, the integrals over $S^0$ and $\partial S^0_f$ appearing in \eqref{23} can be transformed like in \eqref{22} into integrals over $\omega$ and $\partial \omega_f\,$, respectively.

We mention that one can consider more general cases of external loads in the definition of the loading potential \eqref{23}, such as for example tracking loads.

\subsection{Main result: Existence of minimizers}

This theorem states the existence of minimizers to the minimization problem \eqref{20}--\eqref{23}.
\begin{theorem}\label{th1}
Assume that the external loads satisfy the conditions
\begin{equation}\label{24}
    \boldsymbol{f}\in\boldsymbol{L}^2(\omega,\mathbb{R}^3),\qquad  \boldsymbol{n}^*\in \boldsymbol{L}^2(\partial\omega_f,\mathbb{R}^3),
\end{equation}
and the boundary data satisfy the conditions
\begin{equation}\label{25}
    \boldsymbol{y}^*\in\boldsymbol{H}^1(\omega ,\mathbb{R}^3),\qquad \boldsymbol{R}^*\in\boldsymbol{H}^1(\omega, SO(3)).
\end{equation}
Assume that the following conditions concerning the initial configuration are fulfilled: $\,\boldsymbol{y}^0:\omega\subset \mathbb{R}^2\rightarrow\mathbb{R}^3$ is a continuous injective mapping and
\begin{equation}\label{26}
    \begin{array}{c}
   \boldsymbol{y}^0\in\boldsymbol{H}^1(\omega ,\mathbb{R}^3),\qquad \boldsymbol{Q}^{0}\in\boldsymbol{H}^1(\omega, SO(3)),
   \end{array}
\end{equation}
\begin{equation}\label{26,1}
    \begin{array}{c}
   \boldsymbol{a}_\alpha= \partial_\alpha \boldsymbol{y}^0\in \boldsymbol{L}^\infty(\omega ,\mathbb{R}^3)\quad \big(\,\mathrm{i.e.}\,\,\, \nabla\boldsymbol{y}^0\in \boldsymbol{L}^\infty(\omega ,\mathbb{R}^{3\times 2})\big), \vspace{6pt}\\
    \det\big(a_{\alpha\beta}(x_1,x_2)\big)\geq a_0^2 >0\,\,,
    \end{array}
\end{equation}
where $a_0$ is a constant. The strain energy density $W(\boldsymbol{E}^e,\boldsymbol{K}^e)$ is assumed to be a quadratic convex function of $(\boldsymbol{E}^e,\boldsymbol{K}^e)$ and $W$ is coercive, in the sense that  there exists a constant $C_0>0$ with
\begin{equation}\label{26bis}
    W( {\boldsymbol{E}^e}, {\boldsymbol{K}^e})\,\geq\, C_0\,\big( \,   \|\boldsymbol{E}^e\|^2 +  \|\boldsymbol{K}^e\|^2\,\big).
\end{equation}
Then, the minimization problem \eqref{20}--\eqref{23} admits at least one minimizing solution pair
$(\hat{\boldsymbol{y}},\hat{\boldsymbol{R}})\in  \mathcal{A}$.
\end{theorem}
\begin{remark}
The hypotheses \eqref{26,1} can be written equivalently in terms of the tensor $\boldsymbol{P}=\nabla \, \boldsymbol{\Theta}(x_1,x_2,0)$ as
\begin{equation}\label{26,2}
    \begin{array}{c}
   \boldsymbol{P} \in \boldsymbol{L}^\infty(\omega ,\mathbb{R}^{3\times 3}),\qquad \det\boldsymbol{P} \geq a_0 >0\,\,,
    \end{array}
\end{equation}
in view of the relations \eqref{14} and \eqref{14bis}. Since $\boldsymbol{y}^0$ represents the position vector of the reference base surface $S^0$ (which is bounded), the conditions \eqref{26}$_1$ and \eqref{26,1}$_1$ can be written together in the form $\boldsymbol{y}^0 \in\boldsymbol{W}^{1,\infty}(\omega ,\mathbb{R}^3)$.
\end{remark}
\begin{proof}
We employ the direct methods of the calculus of variations. We show first that there exists a constant $C>0$ such that
\begin{equation}\label{27}
    |\,\Lambda(\boldsymbol{y},\boldsymbol{R})\,|\,\leq\, \,\,C\,\big(\,\|\boldsymbol{y}\|_{H^1(\omega)}+1\big),\quad\forall\,(\boldsymbol{y},\boldsymbol{R})\in\mathcal{A}.
\end{equation}
Indeed, since $\boldsymbol{a}_\alpha\in\boldsymbol{L}^\infty(\omega,\mathbb{R}^3)$ it follows that $a=\sqrt{\det(a_{\alpha\beta})}\in{L}^\infty(\omega)$.
We also have $\|\boldsymbol{R}\|^2=\text{tr}(\boldsymbol{R}\boldsymbol{R}^T)=3$, $\,\forall\boldsymbol{R}\in SO(3)$. Taking into account the hypotheses \eqref{24} and the boundedness of $\,\,\Pi_{S^0}\,$ and $ \,\Pi_{\partial S^0_f}\,$, we deduce from \eqref{23} that
\begin{equation*}
 \begin{array}{l}
     |\,\Lambda(\boldsymbol{y},\boldsymbol{R})\,|\,\leq\, |\, \Lambda_{S^0}(\boldsymbol{y},\boldsymbol{R})\, | + |\,
    \Lambda_{\partial S^0_f}(\boldsymbol{y},\boldsymbol{R}) \,|\,\leq   \,
     C_1\,\|\boldsymbol{y}-\boldsymbol{y}^0\|_{L^2(\omega)} +C_2\, \|\boldsymbol{y}-\boldsymbol{y}^0\|_{L^2(\partial\omega_f)} \\
     \qquad +\,|\, \Pi_{S^0}(\boldsymbol{R})\,|+ | \,\Pi_{\partial S^0_f}(\boldsymbol{R})\,| \,\leq \, C_3 \|\boldsymbol{y}\|_{L^2(\omega)}+C_4\|\boldsymbol{y}\|_{H^1(\omega)} +C_5\,,
     \end{array}
\end{equation*}
for some positive constants $C_k>0$. Then, the inequality \eqref{27} holds.

In what follows, we employ the component form of the elastic  strain tensors $\boldsymbol{E}^e$ and $\boldsymbol{K}^e$, written as matrices $E^e$ and $K^e$ in \eqref{16} and \eqref{18}, respectively. Let us show next that there exists a positive constant $\lambda_0>0$ such that
\begin{equation}\label{28}
    \|\,\boldsymbol{E}^e\,\|=\|\,E^e\,\|\,\geq\, \lambda_0\,\|\,H\,\|\,,\qquad  \|\,\boldsymbol{K}^e\,\|=\|\,K^e\,\|\,\geq\, \lambda_0\,\|\,L\,\|\,,
\end{equation}
where the matrices $H=\big(H_{ij}\big)_{3\times 3}$ and $L=\big(L_{ij}\big)_{3\times 3}$ are introduced in \eqref{16} and \eqref{18}. Indeed, since $E^e=Q^0HP^{-1}$ and $Q^0\in SO(3)$ we have
\begin{equation}\label{29}
    \|\,E^e\,\|^2=\|\,Q^0HP^{-1}\,\|^2=\|\,HP^{-1}\,\|^2=\text{tr}\big[ HP^{-1} (HP^{-1})^T \big]=\text{tr}\big[ H(P^{T}\!P)^{-1} H^T \big].
\end{equation}
From \eqref{14bis} we deduce that
\begin{equation}\label{30}
    P^{T}P =\begin{bmatrix} a_{11} & a_{12} &  0 \\
                                 a_{12} & a_{22} &  0 \\
                                 0 & 0 & 1
    \end{bmatrix}\qquad\text{and therefore}\qquad
    (P^{T}P)^{-1}=\begin{bmatrix} a^{11} & a^{12} &  0 \\
                                 a^{12} & a^{22} &  0 \\
                                 0 & 0 & 1
    \end{bmatrix}\,\,.
\end{equation}
Inserting \eqref{30} into \eqref{29} we obtain
\begin{equation}\label{31}
     \|\,E^e\,\|^2=a^{\alpha\beta}H_{i\alpha}H_{i\beta} = a^{\alpha\beta}H_{1\alpha}H_{1\beta}+ a^{\alpha\beta}H_{2\alpha}H_{2\beta}+ a^{\alpha\beta}H_{3\alpha}H_{3\beta}\,,\,\,\,\text{with}\,\,\, H=\big(H_{ij}\big)_{3\times 3}\,,
\end{equation}
since $H_{i3}=0$ according to \eqref{16}.
In virtue of \eqref{26,1}, it follows that the matrix $\big(a_{\alpha\beta}\big)_{2\times 2}$ and its inverse matrix
$\big(a^{\alpha\beta}\big)_{2\times 2}=\big(a_{\alpha\beta}\big)^{-1}$  satisfy
$$\big(a_{\alpha\beta}\big)
\in\boldsymbol{L}^\infty(\omega, \mathbb{R}^{2\times 2})\quad\text{and}\quad \big(a^{\alpha\beta}\big)
\in\boldsymbol{L}^\infty(\omega, \mathbb{R}^{2\times 2}).$$
Then, the smallest eigenvalue of the positive definite symmetric matrix $\big(a^{\alpha\beta}(x_1,x_2)\big)_{2\times 2}$ is greater than a positive constant $\lambda_0^2>0$ and consequently
\begin{equation}\label{32}
    a^{\alpha\beta}\!(x_1,x_2)\,\,v_\alpha\, v_\beta\,\geq\, \lambda_0^2\,v_\gamma\, v_\gamma\,,\qquad \forall\,(x_1,x_2)\in\omega,\,\,\, \forall\,v_1,v_2\in\mathbb{R}.
\end{equation}
Using  inequality \eqref{32} for each individual sum in the right-hand side of \eqref{31} we deduce that $\|E^e\|^2\geq \lambda_0^2\,H_{i\alpha}H_{i\alpha}=\lambda_0^2 \,\|H\|^2$, i.e. the inequality \eqref{28}$_1$ is proved. The proof of the inequality \eqref{28}$_2$ is identical. In view of \eqref{28}$_1$ and \eqref{16} we have
\begin{equation*}
\begin{array}{l}
\|\,E^e\,\|^2\geq \lambda_0^2\displaystyle{\sum_{\alpha=1}^2} \|\,R^T\partial_\alpha y\!-\! Q^{0,T}\partial_\alpha y^{0}\,\|^2   \\
 \qquad\qquad=
 \lambda_0^2\displaystyle{\sum_{\alpha=1}^2}\big( \|\,R^T\partial_\alpha y\|^2\!-2\langle\,   R^T\partial_\alpha y\, , Q^{0,T}\partial_\alpha y^{0}\,\rangle +\|\,Q^{0,T}\partial_\alpha y^{0} \|^2\big)\\
 \qquad\qquad = \lambda_0^2\displaystyle{\sum_{\alpha=1}^2}\big( \|\,\partial_\alpha y\|^2-2\langle\,   R^T\partial_\alpha y\, , Q^{0,T}\partial_\alpha y^{0}\,\rangle +\|\,\partial_\alpha y^{0} \|^2\big)
 ,
 \end{array}
\end{equation*}
where $\langle\,S,T\,\rangle=\text{tr}[ST^T]$ is the scalar product of  two matrices $S,T$. Integrating over $\omega$ and using the Cauchy--Schwarz inequality we obtain
$$\|\,E^e\,\|^2_{L^2(\omega)} \geq
 \lambda_0^2\sum_{\alpha=1}^2\big( \|\,\partial_\alpha y\|^2_{L^2(\omega)}-2\|\, \partial_\alpha y\|_{L^2(\omega)} \|\,\partial_\alpha y^{0}\,\|_{L^2(\omega)} +\|\,\partial_\alpha y^{0 } \|^2_{L^2(\omega)}\big),$$
or
\begin{equation}\label{33}
    \|\,\boldsymbol{E}^e\,\|^2_{L^2(\omega)} \geq
 \lambda_0^2\,\big( \|\,\partial_1 \boldsymbol{y}\|^2_{L^2(\omega)} +\|\,\partial_2 \boldsymbol{y} \|^2_{L^2(\omega)}\big) - \bar{C}_1\|\, \boldsymbol{y}\,\|_{H^1(\omega)}+\bar{C}_2\,,
\end{equation}
for some positive constants $\bar{C}_1>0$, $\bar{C}_2>0$. Let us show that the functional $I(\boldsymbol{y},\boldsymbol{R})$ is bounded from below over the admissible set $\mathcal{A}$. By virtue of \eqref{22}, \eqref{26,1}$_2$ and \eqref{27} we can write
$$I(\boldsymbol{y},\boldsymbol{R})\geq C_0 \!\!\int_\omega \|\,\boldsymbol{E}^e\,\|^2\,a\, \mathrm{d}x_1\mathrm{d}x_2 - \Lambda(\boldsymbol{y},\boldsymbol{R}) \geq C_0\,a_0 \|\,\boldsymbol{E}^e\,\|^2_{L^2(\omega)}  - C\,\big(\,\|\boldsymbol{y}\|_{H^1(\omega)}+1\big)$$
and using \eqref{33} we deduce that there exist the constants $\bar{C}_3>0$ and $\bar{C}_4$ such that
\begin{equation}\label{34}
    I(\boldsymbol{y},\boldsymbol{R})\geq C_0\,a_0 \lambda_0^2\,\big( \|\,\partial_1 \boldsymbol{y}\|^2_{L^2(\omega)} +\|\,\partial_2 \boldsymbol{y} \|^2_{L^2(\omega)}\big) - \bar{C}_3\|\, \boldsymbol{y}\,\|_{H^1(\omega)}-\bar{C}_4\,,\quad\forall\,(\boldsymbol{y},\boldsymbol{R})\in \mathcal{A},
\end{equation}
with $\,a_0\,$ specified by \eqref{26,1}.
We observe that the vector field $\boldsymbol{y}-\boldsymbol{y}^*\in\boldsymbol{H}^1(\omega,\mathbb{R}^3)$ satisfies $\boldsymbol{y}-\boldsymbol{y}^*=\boldsymbol{0}$ on $\partial\omega_d\,$. Applying the Poincar\'e--inequality we infer the existence of a constant $c_p>0$ such that
\begin{equation}\label{35}
     \| \partial_1 (\boldsymbol{y}-\boldsymbol{y}^*)\|^2_{L^2(\omega)} +\| \partial_2 (\boldsymbol{y}-\boldsymbol{y}^*) \|^2_{L^2(\omega)}  \,\geq\, c_p\,\|\, \boldsymbol{y}-\boldsymbol{y}^*\,\|_{H^1(\omega)}^2\,\,.
\end{equation}
Using inequalities of the type $\|\partial_\alpha \boldsymbol{y}\|_{L^2(\omega)}^2\geq \big(\|\partial_\alpha (\boldsymbol{y}-\boldsymbol{y}^*)\|_{L^2(\omega)}- \| \partial_\alpha \boldsymbol{y}^*\|_{L^2(\omega)}\big)^2$ and \eqref{35} we find that
\begin{equation*}
\begin{array}{l}
\|\partial_1 \boldsymbol{y}\|^2_{L^2(\omega)} +\|\partial_2 \boldsymbol{y} \|^2_{L^2(\omega)}\geq  c_p\,\| \boldsymbol{y}-\boldsymbol{y}^*\|_{H^1(\omega)}^2 \\
\qquad\quad\,\,\,\, -2\| \boldsymbol{y}-\boldsymbol{y}^*\|_{H^1(\omega)}\big( \|\partial_1 \boldsymbol{y}^*\|_{L^2(\omega)} +\|\partial_2 \boldsymbol{y}^* \|_{L^2(\omega)}\big)  +
\big(\|\,\partial_1 \boldsymbol{y}^*\|^2_{L^2(\omega)} +\|\,\partial_2 \boldsymbol{y}^* \|^2_{L^2(\omega)}\big).
\end{array}
\end{equation*}
From the last inequality and \eqref{34} follows that there exist some constants $\bar{C}_5>0$ and $\bar{C}_6$ with
\begin{equation}\label{36}
    I(\boldsymbol{y},\boldsymbol{R})\geq C_0\,a_0 \lambda_0^2\,c_p\,\| \boldsymbol{y}-\boldsymbol{y}^*\|_{H^1(\omega)}^2-\bar{C}_5\| \boldsymbol{y}-\boldsymbol{y}^*\|_{H^1(\omega)} +\bar{C}_6\,,\quad\forall\,(\boldsymbol{y},\boldsymbol{R})\in \mathcal{A}.
\end{equation}
Since the constant $C_0\,a_0 \lambda_0^2\,c_p>0$, the function $I(\boldsymbol{y},\boldsymbol{R})$ is bounded from below over   $\mathcal{A}$.  Hence, there exists an infimizing sequence $\big\{(\boldsymbol{y}_n,\boldsymbol{R}_n)\big\}_{n=1}^\infty \subset\mathcal{A}$ such that
\begin{equation}\label{37}
    \lim_{n\rightarrow \infty} I(\boldsymbol{y}_n,\boldsymbol{R}_n) = \,\inf\, \big\{I(\boldsymbol{y},\boldsymbol{R})\, \big|\,  (\boldsymbol{y},\boldsymbol{R})\in \mathcal{A}\big\}.
\end{equation}
In view of the conditions \eqref{25} we have $I(\boldsymbol{y}^*,\boldsymbol{R}^*)<\infty$. The infimizing sequence $\big\{(\boldsymbol{y}_n,\boldsymbol{R}_n)\big\}_{n=1}^\infty$ can be chosen such that
\begin{equation}\label{38}
    I(\boldsymbol{y}_n,\boldsymbol{R}_n)\,\leq \,I(\boldsymbol{y}^*,\boldsymbol{R}^*)\,< \infty\,, \qquad \forall\,n\geq 1.
\end{equation}
Taking into account \eqref{36} and \eqref{38} we see that the sequence $\big\{ \boldsymbol{y}_n \big\}_{n=1}^\infty$ is bounded in $\boldsymbol{H}^1(\omega,\mathbb{R}^3)$. Then, we can extract a subsequence of $\big\{ \boldsymbol{y}_n \big\}_{n=1}^\infty$  (not relabeled) which converges weakly in  $\boldsymbol{H}^1(\omega,\mathbb{R}^3)$ and moreover, according to Rellich's selection principle, it converges strongly in $\boldsymbol{L}^2(\omega,\mathbb{R}^{3})$, i.e. there exists an element $\hat{\boldsymbol{y}}\in\boldsymbol{H}^1(\omega,\mathbb{R}^3)$ such that
\begin{equation}\label{39}
    \boldsymbol{y}_n  \rightharpoonup \hat{ \boldsymbol{y}} \quad\mathrm{in}\quad \boldsymbol{H}^1(\omega, \mathbb{R}^3),\qquad \mathrm{and}\qquad \boldsymbol{y}_n \rightarrow\hat{ \boldsymbol{y}} \quad\mathrm{in}\quad \boldsymbol{L}^2(\omega, \mathbb{R}^3).
\end{equation}
For any $n\in\mathbb{N}$, let us denote by $\boldsymbol{E}^e_n$ and $\boldsymbol{K}^e_n$ the strain measures corresponding to the fields $(\boldsymbol{y}_n,\boldsymbol{R}_n)$, defined by the relations \eqref{8} and \eqref{12}. We have $\boldsymbol{E}^e_n, \boldsymbol{K}^e_n\in \boldsymbol{L}^2(\omega,\mathbb{R}^{3\times 3})$ and let $E_n^e,K_n^e$ be the matrices of components in the basis $\{\boldsymbol{e}_i\otimes\boldsymbol{e}_j\}$, given by \eqref{16} and \eqref{18} in terms $\{y_n,R_n\}$. From \eqref{20}, \eqref{26,1}$_2$ , \eqref{26bis}, \eqref{27} and \eqref{38} we get
$$ C_0a_0\,\|\,\boldsymbol{K}^e_n\,\|^2_{L^2(\omega)} \leq
\int_\omega W(\boldsymbol{E}^e_n,\boldsymbol{K}^e_n)\,a(x_1,x_2)\, \mathrm{d}x_1\mathrm{d}x_2 \leq I(\boldsymbol{y}^*,\boldsymbol{R}^*)+ C\,\big(\,\|\boldsymbol{y}_n\|_{H^1(\omega)}+1\big).
$$
Since $\big\{ \boldsymbol{y}_n \big\}_{n=1}^\infty$ is bounded in $\boldsymbol{H}^1(\omega,\mathbb{R}^3)$, it follows from the last inequalities that $\big\{ \boldsymbol{K}^e_n \big\}_{n=1}^\infty$ is bounded in $\boldsymbol{L}^2(\omega,\mathbb{R}^{3\times 3})$. In view of \eqref{28}$_2\,$, we deduce that $\big\{\text{axl}(\boldsymbol{R}_n^{T}\partial_\alpha \boldsymbol{R}_n)\big\}_{n=1}^\infty$ is bounded in $\boldsymbol{L}^2(\omega,\mathbb{R}^{3})$, or equivalently $ \big\{ \partial_\alpha \boldsymbol{R}_n\big\}_{n=1}^\infty$ is bounded in $\boldsymbol{L}^2(\omega,\mathbb{R}^{3\times 3})$, for $\alpha=1,2$. Since $\boldsymbol{R}_n\in SO(3)$ we have $\|\boldsymbol{R}_n\|^2=3$ and thus we can infer that the sequence $ \big\{ \boldsymbol{R}_n\big\}_{n=1}^\infty$
is bounded in $\boldsymbol{H}^1(\omega,\mathbb{R}^{3\times 3})$. Hence, there exists a subsequence of $ \big\{ \boldsymbol{R}_n\big\}_{n=1}^\infty$ (not relabeled) and an element $\hat{\boldsymbol{R}}\in \boldsymbol{H}^1(\omega,\mathbb{R}^{3\times 3})$ with
\begin{equation}\label{40}
    \boldsymbol{R}_n \rightharpoonup       \hat{\boldsymbol{R}}    \quad\mathrm{in}\quad \boldsymbol{H}^1(\omega, \mathbb{R}^{3\times3}) , \qquad\mathrm{and}\qquad          \boldsymbol{R}_n  \rightarrow    \hat{\boldsymbol{R}} \quad\mathrm{in}\quad \boldsymbol{L}^2(\omega, \mathbb{R}^{3\times3}).
\end{equation}
We can show for the limit that $\hat{\boldsymbol{R}}\in SO(3)$. Indeed, since $\boldsymbol{R}_n\in SO(3)$ we have
$$\|\,\boldsymbol{R}_n \hat{\boldsymbol{R}}{}^T -\id_3\|_{L^2(\omega)}= \|\,\boldsymbol{R}_n ( \hat{\boldsymbol{R}}{}^T -\boldsymbol{R}_n^T)\|_{L^2(\omega)}=  \|\,  \hat{\boldsymbol{R}}  -\boldsymbol{R}_{n }\|_{L^2(\omega)} \longrightarrow 0,
$$
i.e. $\boldsymbol{R}_n\hat{\boldsymbol{R}}{}^T\rightarrow\id_3$ in $\boldsymbol{L}^2(\omega,\mathbb{R}^{3\times 3})$. On the other hand, we can write
$$\|\,\boldsymbol{R}_n \hat{\boldsymbol{R}}{}^T -\hat{\boldsymbol{R}} \hat{\boldsymbol{R}}{}^T\|_{L^1(\omega)}= \|\, ( \boldsymbol{R}_{n }  -\hat{\boldsymbol{R}})\hat{\boldsymbol{R}}{}^T\|_{L^1(\omega)} \leq 3  \|\,   \boldsymbol{R}_{n }  -\hat{\boldsymbol{R}} \|_{L^2(\omega)}\,  \|\,  \hat{\boldsymbol{R}} \|_{L^2(\omega)} \longrightarrow 0,
$$
which means that $\boldsymbol{R}_n\hat{\boldsymbol{R}}{}^T\rightarrow\hat{\boldsymbol{R}}\hat{\boldsymbol{R}}{}^T$ in $\boldsymbol{L}^1(\omega,\mathbb{R}^{3\times 3})$. Consequently, we find $\hat{\boldsymbol{R}}\hat{\boldsymbol{R}}{}^T=\id_3$ so that $\hat{\boldsymbol{R}}\in\boldsymbol{H}^1(\omega,SO(3))$.

By virtue of the relations  $(\boldsymbol{y}_n,\boldsymbol{R}_n)\in \mathcal{A}$ and \eqref{39}, \eqref{40},   we derive that $\hat{\boldsymbol{y}}=\boldsymbol{y}^*$ on $\partial S^0_d$ and $\hat{\boldsymbol{R}}=\boldsymbol{R}^*$ on $\partial S^0_d\,$ in the sense of traces. Hence, we obtain that the limit pair satisfies $(\hat{\boldsymbol{y}},\hat{\boldsymbol{R}})\in \mathcal{A}$.

Let us construct the elements $\hat{\boldsymbol{E}^e}, \hat{\boldsymbol{K}^e} \in \boldsymbol{L}^2(\omega,\mathbb{R}^{3\times 3})$ defined in terms of the fields $(\hat{\boldsymbol{y}},\hat{\boldsymbol{R}})$ by the relations \eqref{8} and \eqref{12}. Then, the matrices of components $\hat{E}^e,\hat{K}^e$ are expressed in terms of the components $(\hat{y},\hat{R})$ by \eqref{16} and \eqref{18}, i.e.
\begin{equation}\label{41}
    \begin{array}{l}
     \hat{E}^e= Q^0\big(\,\hat{R}{}^T\partial_1 \hat{y}- Q^{0,T}\partial_1 y^{0}\,\big|\, \hat{R}{}^T\partial_2 \hat{y}- Q^{0,T}\partial_2 y^{0}\,\big|\,\, 0\,\,\big) \,P^{-1} ,\vspace{4pt}\\
     \hat{K}^e=Q^0\big(\,\text{axl}(\hat{R}{}^T\partial_1 \hat{R})\!-\!\text{axl}(Q^{0,T}\partial_1 Q^{0}) \,\,\big|\,\, \text{axl}(\hat{R}{}^T\partial_2 \hat{R})\! -\!\text{axl}(Q^{0,T}\partial_2 Q^{0})\,\,\big|\, \,0\,\big)P^{-1} \!.
     \end{array}
\end{equation}

Next, we want to show that there exist some subsequences (not relabeled) of $\{\boldsymbol{E}^e_n\}$ and $\{\boldsymbol{K}^e_n\}$ such that
\begin{equation}\label{42}
    \boldsymbol{E}^e_n  \rightharpoonup  \hat{\boldsymbol{E}^e} \quad\mathrm{in}\quad \boldsymbol{L}^2(\omega, \mathbb{R}^{3\times3}),\quad\text{and} \qquad  \boldsymbol{K}^e_n \rightharpoonup \hat{\boldsymbol{K}^e}\quad\mathrm{in}\quad \boldsymbol{L}^2(\omega, \mathbb{R}^{3\times3}).
\end{equation}
As shown above, the sequence $ \big\{ \boldsymbol{y}_n\big\}_{n=1}^\infty$ is bounded in $\boldsymbol{H}^1(\omega,\mathbb{R}^{3 })$. It follows that $ \big\{ \partial_\alpha \boldsymbol{y}_n\big\}_{n=1}^\infty$ is bounded in $\boldsymbol{L}^2(\omega,\mathbb{R}^{3 })$ and then the sequence $ \big\{ \boldsymbol{R}_n^T\partial_\alpha \boldsymbol{y}_n\big\}_{n=1}^\infty$ is bounded in $\boldsymbol{L}^2(\omega,\mathbb{R}^{3 })$, since $ \boldsymbol{R}_n\in SO(3)$. Consequently, there exists a subsequence (not relabeled) and an element $\boldsymbol{\xi}_\alpha\in \boldsymbol{L}^2(\omega,\mathbb{R}^{3})$ such that
\begin{equation}\label{43}
    \boldsymbol{R}_n^T\partial_\alpha \boldsymbol{y}_n \,\,\rightharpoonup \,\, \boldsymbol{\xi}_\alpha\qquad\text{in} \quad \boldsymbol{L}^2(\omega,\mathbb{R}^{3}).
\end{equation}
On the other hand, let $\boldsymbol{\phi}\in  \boldsymbol{C}_0^\infty(\omega,\mathbb{R}^{3})$ be an arbitrary test function. Then, using the properties of the scalar product  we deduce
\begin{equation*}
    \begin{array}{l}
    \displaystyle{\int_\omega}\big(\boldsymbol{R}_n^T \partial_\alpha \boldsymbol{y}_n - \hat{\boldsymbol{R}}{}^T \partial_\alpha \hat{\boldsymbol{y}}\big)\cdot\boldsymbol{\phi} \,\mathrm{d}x_1\mathrm{d}x_2
    \vspace{4pt}\\
    \qquad\quad=
    \displaystyle{\int_\omega}  \hat{\boldsymbol{R}}{}^T \big( \partial_\alpha \boldsymbol{y}_n - \partial_\alpha \hat{\boldsymbol{y}}\big) \cdot\boldsymbol{\phi} \,\mathrm{d}x_1\mathrm{d}x_2
        +\displaystyle{\int_\omega}  \big( \boldsymbol{R}_n^T- \hat{\boldsymbol{R}}{}^T\big) \partial_\alpha \boldsymbol{y}_n  \cdot\boldsymbol{\phi} \,\mathrm{d}x_1\mathrm{d}x_2
    \vspace{4pt}\\
     \qquad\quad   =\displaystyle{\int_\omega}  \big( \partial_\alpha \boldsymbol{y}_n - \partial_\alpha \hat{\boldsymbol{y}}\big) \cdot\hat{\boldsymbol{R}}\boldsymbol{\phi} \, \mathrm{d}x_1\mathrm{d}x_2  +
        \displaystyle{\int_\omega} \!\!\big\langle \boldsymbol{R}_{n}\!\!-\! \hat{\boldsymbol{R}}\,, \partial_\alpha \boldsymbol{y}_n \!\otimes\!\boldsymbol{\phi}\rangle \mathrm{d}x_1\mathrm{d}x_2
    \vspace{4pt}\\
    \qquad\qquad\,\, \leq\|\boldsymbol{R}_{n}\!-\! \hat{\boldsymbol{R}}\|_{L^2(\omega)}\|\partial_\alpha \boldsymbol{y}_n \!\otimes\!\boldsymbol{\phi}\|_{L^2(\omega)} \!+\!\! \displaystyle{\int_\omega}  \big( \partial_\alpha \boldsymbol{y}_n \!- \!\partial_\alpha \hat{\boldsymbol{y}}\big) \!\cdot\!\hat{\boldsymbol{R}}\boldsymbol{\phi} \, \mathrm{d}x_1\mathrm{d}x_2 \,,
\end{array}
\end{equation*}
since the relations \eqref{39}, \eqref{40} and $\hat{\boldsymbol{R}}\boldsymbol{\phi}\in\boldsymbol{L}^2(\omega,\mathbb{R}^{3})$ hold, and $\|\partial_\alpha \boldsymbol{y}_n\otimes\boldsymbol{\phi}\|$ is bounded. Thus, we get
\begin{equation}\label{44}
    \displaystyle{\int_\omega}\big(\boldsymbol{R}_n^T \partial_\alpha \boldsymbol{y}_n \big)\cdot\boldsymbol{\phi} \,\mathrm{d}x_1\mathrm{d}x_2 \longrightarrow
    \displaystyle{\int_\omega}\big( \hat{\boldsymbol{R}}{}^T \partial_\alpha \hat{\boldsymbol{y}}\big)\cdot\boldsymbol{\phi} \,\mathrm{d}x_1\mathrm{d}x_2,\quad\forall\, \boldsymbol{\phi}\in  \boldsymbol{C}_0^\infty(\omega,\mathbb{R}^3).
\end{equation}
By comparison of \eqref{43} and \eqref{44} we find $\boldsymbol{\ell}_\alpha= \hat{\boldsymbol{R}}{}^{T}\partial_\alpha \hat{\boldsymbol{y}}\,$, which means that
$\boldsymbol{R}_n^T \partial_\alpha \boldsymbol{y}_n   \rightharpoonup
     \hat{\boldsymbol{R}}{}^T \partial_\alpha \hat{\boldsymbol{y}}$ in $ \boldsymbol{L}^2(\omega,\mathbb{R}^3)$ ,
or equivalently
\begin{equation}\label{45}
    \big(\boldsymbol{R}_n^T \partial_\alpha \boldsymbol{y}_n  - \boldsymbol{Q}^{0,T} \partial_\alpha \boldsymbol{y}^0 \big)\quad \rightharpoonup\quad
    \big( \hat{\boldsymbol{R}}{}^T \partial_\alpha \hat{\boldsymbol{y}} -  {\boldsymbol{R}}{}^T_0 \partial_\alpha \boldsymbol{y}^{0}\big)\quad\mathrm{in}\quad \boldsymbol{L}^2(\omega,\mathbb{R}^3).
\end{equation}
Taking into account \eqref{16}, \eqref{41}$_1$ and the hypotheses \eqref{26}, \eqref{26,1}, we obtain from \eqref{45} that $E^e_n\rightharpoonup\hat{E}^e$ in $\boldsymbol{L}^2(\omega,\mathbb{R}^{3\times 3})$, i.e. the relation \eqref{42}$_1$ holds.

To prove \eqref{42}$_2$ we start from the fact that the sequence $ \big\{ \boldsymbol{R}_n^T\partial_\alpha \boldsymbol{R}_n\big\}_{n=1}^\infty$ is bounded in $\boldsymbol{L}^2(\omega,\mathbb{R}^{3\times 3 })$, as we proved previously. Then, there exists a subsequence (not relabeled) and an element $\boldsymbol{ \zeta}_\alpha\in \boldsymbol{L}^2(\omega,\mathbb{R}^{3\times 3})$ such that
\begin{equation}\label{46}
    \boldsymbol{R}_n^T\partial_\alpha \boldsymbol{R}_n \rightharpoonup \boldsymbol{ \zeta}_\alpha\qquad\text{in} \quad \boldsymbol{L}^2(\omega,\mathbb{R}^{3\times 3}).
\end{equation}
On the other hand, for any test function $\boldsymbol{\mathit{\Phi}}\in  \boldsymbol{C}_0^\infty(\omega,\mathbb{R}^{3\times 3})$ we can write
\begin{equation*}
\begin{array}{l}
    \displaystyle{\int_\omega}\big\langle \boldsymbol{R}_n^T \partial_\alpha \boldsymbol{R}_{n} - \hat{\boldsymbol{R}}{}^T \partial_\alpha \hat{\boldsymbol{R}}  \,,\, \boldsymbol{\mathit{\Phi}} \big\rangle \,\mathrm{d}x_1\mathrm{d}x_2=
    \displaystyle{\int_\omega}\big\langle   \hat{\boldsymbol{R}}{}^T \big(\partial_\alpha  \boldsymbol{R}_{n} - \partial_\alpha\hat{\boldsymbol{R}}\big) \,,\, \boldsymbol{\mathit{\Phi}} \,\big\rangle \,\mathrm{d}x_1\mathrm{d}x_2 \vspace{2pt}\\
     \qquad+
    \displaystyle{\int_\omega}\big\langle \big( \boldsymbol{R}_n^T- \hat{\boldsymbol{R}}{}^T\big) \partial_\alpha \boldsymbol{R}_{n}  \,,\, \boldsymbol{\mathit{\Phi}} \,\big\rangle \,\mathrm{d}x_1\mathrm{d}x_2
    \leq
     \displaystyle{\int_\omega}
     \big\langle \partial_\alpha \boldsymbol{R}_{n} - \partial_\alpha \hat{\boldsymbol{R}}  \,,\, \hat{\boldsymbol{R}}\boldsymbol{\mathit{\Phi}} \,\big\rangle \,\mathrm{d}x_1\mathrm{d}x_2
     \vspace{2pt}\\
     \qquad+
     \|\boldsymbol{R}_{n}- \hat{\boldsymbol{R}}\|_{L^2(\omega)}\,\| \partial_\alpha \boldsymbol{R}_n \boldsymbol{\mathit{\Phi}}^T \, \|_{L^2(\omega)} \longrightarrow 0,
\end{array}
\end{equation*}
since $\hat{\boldsymbol{R}}\boldsymbol{\mathit{\Phi}}\in\boldsymbol{L}^2(\omega,\mathbb{R}^{3\times 3})$, $\|\partial_\alpha \boldsymbol{R}_n \boldsymbol{\mathit{\Phi}}^T\|$ is bounded, and relations \eqref{40} hold. Consequently, we have
$$\displaystyle{\int_\omega}\big\langle \boldsymbol{R}_n^T \partial_\alpha \boldsymbol{R}_{n}\,,\, \boldsymbol{\mathit{\Phi}} \big\rangle \,\mathrm{d}x_1\mathrm{d}x_2 \longrightarrow
    \displaystyle{\int_\omega}\big\langle  \hat{\boldsymbol{R}}{}^T \partial_\alpha \hat{\boldsymbol{R}}\,,\, \boldsymbol{\mathit{\Phi}} \big\rangle \,\mathrm{d}x_1\mathrm{d}x_2,\quad\forall \,\boldsymbol{\mathit{\Phi}}\in  \boldsymbol{C}_0^\infty(\omega,\mathbb{R}^{3\times3}), $$
and by comparison with \eqref{46} we deduce that $\boldsymbol{ \zeta}_\alpha= \hat{\boldsymbol{R}}{}^{T}\partial_\alpha \hat{\boldsymbol{R}}\,$, i.e. the convergence $ \boldsymbol{R}_n^T \partial_\alpha \boldsymbol{R}_n\rightharpoonup \hat{\boldsymbol{R}}{}^{T}\partial_\alpha \hat{\boldsymbol{R}}$ holds in $\boldsymbol{L}^2(\omega,\mathbb{R}^{3\times 3})$. It follows that
$$\big[\text{axl}(\boldsymbol{R}_n^T \partial_\alpha \boldsymbol{R}_n) - \text{axl}(\boldsymbol{R}^{T}_0 \partial_\alpha \boldsymbol{Q}^{0})\big]
 \,\,\,\rightharpoonup \,\,\,
\big[\text{axl}(\hat{\boldsymbol{R}}{}^{T}\partial_\alpha \hat{\boldsymbol{R}})- \text{axl}(\boldsymbol{R}^{T}_0 \partial_\alpha \boldsymbol{Q}^{0})\big]
\quad\text{in}\,\, \boldsymbol{L}^2(\omega,\mathbb{R}^{3\times 3}),$$
and from \eqref{18}, \eqref{26}, \eqref{26,1} and \eqref{41}$_2$ we derive that the convergence \eqref{42}$_2$ holds true. \medskip

In the last step of the proof we use the convexity of the strain energy density $W$. In view of \eqref{42}, we have
\begin{equation}\label{47}
     \int_\omega W(\hat{\boldsymbol{E}^e},\hat{\boldsymbol{K}^e})\,a(x_1,x_2)\,\mathrm{d}x_1\mathrm{d}x_2\,\leq \, \liminf_{n\to\infty}  \int_\omega W( {\boldsymbol{E}^e_n}, {\boldsymbol{K}^e_n})\,a(x_1,x_2)\,\mathrm{d}x_1\mathrm{d}x_2.
\end{equation}
since $W$ is convex in $(\boldsymbol{E}^e,\boldsymbol{K}^e)$. Taking into account the hypotheses \eqref{24}, the continuity of the load potential functions      $\,\,\Pi_{S^0}\,$, $ \,\Pi_{\partial S^0_f}\,$,    and the convergence relations \eqref{39}$_2$ and \eqref{40}$_2\,$, we deduce
\begin{equation}\label{48}
  \Lambda(\hat{\boldsymbol{y}}, \hat{\boldsymbol{R}})=  \lim_{n\to\infty} \Lambda(\boldsymbol{y}_n, \boldsymbol{R}_n).
\end{equation}
From \eqref{20}, \eqref{22}, \eqref{47} and \eqref{48} we get
\begin{equation}\label{49}
    I(\hat{\boldsymbol{y}},\hat{\boldsymbol{R}})\,\leq\,  \liminf_{n\to\infty} \, I(\boldsymbol{y}_n,\boldsymbol{R}_n)\,.
\end{equation}
Finally, the relations \eqref{37} and \eqref{49} show that
$$ I(\hat{\boldsymbol{y}},\hat{\boldsymbol{R}})\,=\,
\,\inf\, \big\{I(\boldsymbol{y},\boldsymbol{R})\, \big|\,  (\boldsymbol{y},\boldsymbol{R})\in \mathcal{A}\big\}.
$$
Since $(\hat{\boldsymbol{y}},\hat{\boldsymbol{R}})\in\mathcal{A}$, we conclude that $(\hat{\boldsymbol{y}},\hat{\boldsymbol{R}})$ is a minimizing solution pair of our minimization problem. The proof is complete.
\hfill\end{proof}
\begin{remark}
The solution fields satisfy the following regularity conditions
$$\hat{\boldsymbol{y}}\in \boldsymbol{H}^1(\omega,\mathbb{R}^{3})
,\qquad
\hat{\boldsymbol{R}}\in\boldsymbol{L}^\infty(\omega,SO(3))
\cap \boldsymbol{H}^1(\omega,SO(3)) .
$$
Thus, the position vector $ \,\hat{\boldsymbol{y}}\,$ and the total rotation field $\,\hat{\boldsymbol{R}}\,$ may fail to be continuous, according to the limit case of Sobolev embedding.
\end{remark}
\begin{remark}
We observe that the boundary conditions imposed on the orthogonal tensor $\boldsymbol{R}$ can be relaxed in the definition of the admissible set $\mathcal{A}$. Thus, one can prove the existence of minimizers for the minimization problem \eqref{20} over the following larger admissible set
\begin{equation*}
    \tilde{\mathcal{A}} =\big\{(\boldsymbol{y},\boldsymbol{R})\in\boldsymbol{H}^1(\omega,\mathbb{R}^3) \times\boldsymbol{H}^1(\omega, SO(3))\,\,\, \big|\,\, \,\, \boldsymbol{y}_{\big| \partial\omega_d}=\boldsymbol{y}^*  \big\}.
\end{equation*}
This assertion can be proved in the same way as the Theorem \ref{th1}. For a discussion of   possible alternative boundary conditions for the   field $\boldsymbol{R}$ on $\partial\omega_d$ we refer to the works \cite{Neff_plate04_cmt,Neff_plate07_m3as}.
\end{remark}

\section{Applications of the theorem and discussions}\label{sect5}

In this section we present some important special cases for the choice of the energy density $W$ where the Theorem \ref{th1} can be successfully applied to show the existence of minimizers.

Let us discuss first on the choice of the 3 initial directors $\{\boldsymbol{d}_i^0\}$ in the reference configuration, i.e. the specification of the proper orthogonal tensor $\boldsymbol{Q}^{0}=\boldsymbol{d}_i^0 \otimes \boldsymbol{e}_i\,$. One judicious choice for the tensor $\boldsymbol{Q}^{0}$ is the following
\begin{equation}\label{49,1}
    \boldsymbol{Q}^{0}=\text{polar}(\boldsymbol{P})=\text{polar}\big( \nabla \boldsymbol{\Theta}(x_1,x_2,0)\big),
\end{equation}
where $\boldsymbol{P}=\boldsymbol{a}_{i}\otimes\boldsymbol{e}_i=\partial_\alpha \boldsymbol{y}^0\otimes\boldsymbol{e}_\alpha +\boldsymbol{n}^0\otimes \boldsymbol{e}_3$ has been introduced previously in \eqref{14} and $\text{polar}(\boldsymbol{T})$ denotes the orthogonal tensor given by the polar decomposition of any tensor $\boldsymbol{T}$.
\begin{remark}\label{rem6}
If the tensor $\boldsymbol{Q}^{0}$ is given by \eqref{49,1}, then the (initial) directors $\boldsymbol{d}_\alpha^0$ belong to the tangent plane at any point of $S^0$ and $\boldsymbol{d}_3^0= \boldsymbol{n}^0$. Indeed, let $\boldsymbol{P}=\boldsymbol{Q}^{0}\boldsymbol{U}^0$ be the polar decomposition of $\boldsymbol{P}$.
Using the matrices of components in the $\{\boldsymbol{e}_{i}\otimes\boldsymbol{e}_j\}$ tensor basis, we write this relation as $P=Q^0U^0$, and from \eqref{30} we derive consecutively
\begin{equation*}
    U^{0,T} U^{0}=P^{T}P =\begin{bmatrix} a_{11} & a_{12} &  0 \\
                                 a_{12} & a_{22} &  0 \\
                                 0 & 0 & 1
    \end{bmatrix},\quad
    U^0=
    \begin{bmatrix} u_{11}^0 & u_{12}^0 &  0 \\
                                 u_{12}^0 & u_{22}^0 &  0 \\
                                 0 & 0 & 1
    \end{bmatrix},\quad
    \big(U^0\big)^{-1}=
    \begin{bmatrix} \bar{u}_{11}^0 & \bar{u}_{12}^0 &  0 \\
                                 \bar{u}_{12}^0 & \bar{u}_{22}^0 &  0 \\
                                 0 & 0 & 1
    \end{bmatrix},
\end{equation*}
where $u_{\alpha\beta}^0$ and $\bar{u}_{\alpha\beta}^0$ are some given real functions of $(x_1,x_2)$. In view of \eqref{14bis}, it follows
\begin{equation}\label{49,2}
   Q^0= P\,\big(U^0\big)^{-1}= \Big(\,a_{1}\,\Big|\, a_{2}\,\Big|\, n^0\, \Big)_{3\times 3}
   \begin{bmatrix} \bar{u}_{11}^0 & \bar{u}_{12}^0 &  0 \\
                                 \bar{u}_{12}^0 & \bar{u}_{22}^0 &  0 \\
                                 0 & 0 & 1
    \end{bmatrix}\qquad\Rightarrow\qquad Q^0e_3=n^0,
\end{equation}
from which we can see that the third column of the matrix $\,Q^0$ is equal to $\,n^0$.
On the other hand, by the definition \eqref{1}$_{2\,}$, the initial rotation field $ \boldsymbol{Q}^{0}$ is given by $ \boldsymbol{Q}^{0}= \boldsymbol{d}_i^{0}\otimes  \boldsymbol{e}_i\,$ and the matrix $Q^0$ can be written in column form as
\begin{equation}\label{49,3}
   Q^0 = \Big(\,d_{1}^0\,\Big|\, d_{2}^0\,\Big|\, d_3^0\, \Big)_{3\times 3}\,\,\,.
\end{equation}
If we compare \eqref{49,2} and  \eqref{49,3} we find that $d_3^0=n^0$. Thus, we have
$\boldsymbol{d}_3^{0}= \boldsymbol{n}^{0}$ and  $\{\boldsymbol{d}_1^0,\boldsymbol{d}_2^0\}$ is an orthonormal basis in the tangent plane, at any point of $S^0$.

If we choose the tensor $\boldsymbol{Q}^{0}$ as in \eqref{49,1}, then in order to satisfy \eqref{26}$_2$ we need to consider an additional regularity assumption on the initial configuration, namely
$$\mathrm{polar}(\boldsymbol{P})=\mathrm{polar}\big( \nabla \boldsymbol{\Theta}(x_1,x_2,0)\big)\in \boldsymbol{H}^1(\omega,SO(3)),$$
which is equivalent to $\mathrm{Curl}\big[\mathrm{polar}\big( \nabla \boldsymbol{\Theta}(x_1,x_2,0)\big)\big]\in \boldsymbol{L}^2(\omega,SO(3))\,$, cf. \cite{Neff_curl08}. A stronger sufficient condition is $\boldsymbol{\Theta}\in \boldsymbol{W}^{1,\infty}(\omega,\mathbb{R}^3)\cap  \boldsymbol{H}^2(\omega,\mathbb{R}^3)$.
\end{remark}

It is possible to simplify the form of the equations in the case of an orthogonal parametrization of the initial surface $S^0\,$. If we   assume that the curvilinear coordinates $(x_1,x_2)$ are   such that the basis $\{\boldsymbol{a}_1,\boldsymbol{a}_2,\boldsymbol{n}^0\}$ is orthonormal, then the initial surface $S^0$ is formally parametrized by orthogonal arc-length coordinates \cite{Chroscielewski11} and we have
\begin{equation}\label{50}
    \boldsymbol{a}_\alpha=\boldsymbol{a}^\alpha,\qquad a_{\alpha\beta}= a^{\alpha\beta}=\delta_{\alpha\beta}\,.
\end{equation}
\begin{remark}
The \emph{Theorema Egregium} (Gauss) can be put into the following form: the Gaussian curvature $K$ can be found given the full knowledge of the first fundamental form of the surface and expressed via the first fundamental form and its partial derivatives of first and second order (the Brioschi formula). Therefore, the Gaussian curvature of an embedded smooth surface in $\mathbb{R}^3$ is invariant under local isometries, i.e. if the parametrization $\boldsymbol{y}^0:\omega\subset\mathbb{R}^2\rightarrow \mathbb{R}^3$ of the surface from a flat reference configuration $\omega$ is given such that  $\big(\nabla \boldsymbol{y}^0\big)^T \nabla \boldsymbol{y}^0=\id_2$ (the basis $\{\boldsymbol{a}_1,\boldsymbol{a}_2,\boldsymbol{n}^0\}$ is orthonormal), then the curvature $K$ of the surface $\boldsymbol{y}^0(\omega)$ is necessarily zero. This is only the case for developable surfaces.

For general surfaces it is therefore impossible to determine, even locally, an orthonormal parametrization. However, in FEM approaches one may think in a discrete pointwise manner as in \cite{Chroscielewski11}.

For example, let $S^0$ be a cylindrical surface (which is a developable surface) with generators parallel to $\boldsymbol{e}_3\,$. The position vector $\boldsymbol{y}^0$ is given by
$$\boldsymbol{y}^0=\boldsymbol{y}^0(\theta,z)=r_0\cos\dfrac{\theta}{r_0}\, \boldsymbol{e}_1+
r_0\sin\dfrac{\theta}{r_0}\, \boldsymbol{e}_2+ z \, \boldsymbol{e}_3\qquad (r_0>0\,\,\mathrm{constant}).
$$
Choosing the curvilinear coordinates $x_1=\theta, x_2=z$, we have
$$\boldsymbol{a}_1= \partial_1 \boldsymbol{y}^0=- \sin\dfrac{\theta}{r_0}\, \boldsymbol{e}_1+
\cos\dfrac{\theta}{r_0}\, \boldsymbol{e}_2,\quad
\boldsymbol{a}_2=\partial_2 \boldsymbol{y}^0=\boldsymbol{e}_3
,\quad
\boldsymbol{n}^0=\cos\dfrac{\theta}{r_0}\, \boldsymbol{e}_1+
\sin\dfrac{\theta}{r_0}\, \boldsymbol{e}_2\,,$$
so that $\{\boldsymbol{a}_1,\boldsymbol{a}_2,\boldsymbol{n}^0\}$ is orthonormal.
\end{remark}\smallskip

In view of \eqref{49,1}--\eqref{50}, we obtain in this case that $\boldsymbol{Q}^{0}=\boldsymbol{P}$ (since $\boldsymbol{U}^0=\id_3$ and $\text{polar}(\boldsymbol{P})=  \boldsymbol{P}\in SO(3)$) and the directors $\{\boldsymbol{d}_i^0\}$ in the reference configuration coincide with $\{\boldsymbol{a}_1,\boldsymbol{a}_2,\boldsymbol{n}^0\}$ in each point of $S^0\,$, i.e.
\begin{equation}\label{51}
    \boldsymbol{d}_\alpha^0=\boldsymbol{a}_\alpha= \partial_\alpha \boldsymbol{y}^0\,,\qquad
    \boldsymbol{d}_3^0= \boldsymbol{a}_3=\boldsymbol{n}^0,\qquad \boldsymbol{a}_i= \boldsymbol{Q}^{0}\boldsymbol{e}_i\,,\qquad \boldsymbol{d}_i=\boldsymbol{R} \boldsymbol{Q}^{0,T}\boldsymbol{a}_i\,.
\end{equation}
The expressions of the elastic strain measures $\boldsymbol{E}^e$ and $\boldsymbol{K}^e$ may be simplified  in this situation. By virtue of \eqref{50} and \eqref{51} we get
\begin{equation}\label{52}
\begin{array}{l}
    \boldsymbol{Q}^{0} \big(\boldsymbol{R}^T\partial_\alpha \boldsymbol{y} -  \boldsymbol{Q}^{0,T} \partial_\alpha \boldsymbol{y}^0\big) = \boldsymbol{Q}^{0}  \boldsymbol{R}^T\partial_\alpha \boldsymbol{y} -  \boldsymbol{a}_{ \alpha} \vspace{3pt}\\
    \qquad\qquad\quad\qquad\qquad\qquad
     = \big( \boldsymbol{a}_i\cdot\boldsymbol{Q}^{0}  \boldsymbol{R}^T\partial_\alpha \boldsymbol{y} - \delta_{\alpha i} \big) \boldsymbol{a}_{ i}
     = \big( \boldsymbol{d}_i\cdot\partial_\alpha \boldsymbol{y} - \delta_{\alpha i} \big) \boldsymbol{a}_{ i}\,.
     \end{array}
\end{equation}
We can write
$\boldsymbol{R}^T\partial_\alpha \boldsymbol{R}=(\boldsymbol{d}_{i}\otimes\boldsymbol{e}_i)^T (\partial_\alpha \boldsymbol{d}_{j}\otimes\boldsymbol{e}_j)= (\boldsymbol{d}_{i}\cdot \partial_\alpha \boldsymbol{d}_{j} )
\boldsymbol{e}_i\otimes \boldsymbol{e}_j\,$,
so that we find
$$\text{axl}(\boldsymbol{R}^T\partial_\alpha \boldsymbol{R})= \dfrac{1}{2}\,e_{ijk}(\boldsymbol{d}_{k}\cdot \partial_\alpha \boldsymbol{d}_{j} )\boldsymbol{e}_i\,,\qquad \text{axl}(\boldsymbol{Q}^{0,T} \partial_\alpha \boldsymbol{Q}^{0})= \dfrac{1}{2}\,e_{ijk}(\boldsymbol{a}_{k}\cdot \partial_\alpha \boldsymbol{a}_{j} )\boldsymbol{e}_i\,,
$$
where $e_{ijk}$ is the permutation symbol. The last relations and $\boldsymbol{Q}^{0,T} \boldsymbol{a}_i =\boldsymbol{e}_i$ yield
\begin{equation}\label{53}
\begin{array}{c}
    \boldsymbol{Q}^{0} \big[\text{axl}(\boldsymbol{R}^T\partial_\alpha \boldsymbol{R}) - \text{axl}(\boldsymbol{Q}^{0,T} \partial_\alpha \boldsymbol{Q}^{0})\big]= \big[\big(\text{axl}(\boldsymbol{R}^T\partial_\alpha \boldsymbol{R}) - \text{axl}(\boldsymbol{Q}^{0,T} \partial_\alpha \boldsymbol{Q}^{0})\big)\!\cdot\! \boldsymbol{e}_i\big]\boldsymbol{a}_i \\
    =
    \dfrac{1}{2}\,e_{ijk}\big[ (\boldsymbol{d}_{k}\cdot \partial_\alpha \boldsymbol{d}_{j} ) - (\boldsymbol{a}_{k}\cdot \partial_\alpha \boldsymbol{a}_{j} )\big]\boldsymbol{a}_i\,.
    \end{array}
\end{equation}
Using the relations \eqref{8}, \eqref{12}, \eqref{52} and \eqref{53} we decompose the strain tensor $\boldsymbol{E}^e$ and the curvature tensor $\boldsymbol{K}^e$ in the basis $\{\boldsymbol{a}_i\otimes \boldsymbol{a}_\alpha\}$ as follows
\begin{equation}\label{54}
\begin{array}{l}
    \boldsymbol{E}^e=\tilde{E}^e_{i\alpha}\boldsymbol{a}_i\otimes\boldsymbol{a}_\alpha\,,\qquad
    \tilde{E}^e_{i\alpha}= \boldsymbol{d}_i  \cdot   \partial_\alpha \boldsymbol{y} - \delta_{\alpha i}\,,\qquad
        \boldsymbol{K}^e=\tilde{K}^e_{i\alpha}\boldsymbol{a}_i\otimes\boldsymbol{a}_\alpha\,,\\
        \tilde{K}^e_{1\alpha}= \boldsymbol{d}_3\!\cdot\!\partial_\alpha \boldsymbol{d}_{2} \! -\! \boldsymbol{a}_3\!\cdot\!\partial_\alpha \boldsymbol{a}_{2}, \quad
    \tilde{K}^e_{2\alpha}= \boldsymbol{d}_1\!\cdot\!\partial_\alpha \boldsymbol{d}_{3} \!- \! \boldsymbol{a}_1\!\cdot\!\partial_\alpha \boldsymbol{a}_{3}, \quad\!
    \tilde{K}^e_{3\alpha}= \boldsymbol{d}_2\!\cdot\!\partial_\alpha \boldsymbol{d}_{1} \!-\! \boldsymbol{a}_2 \!\cdot\!\partial_\alpha \boldsymbol{a}_{1} .
\end{array}
\end{equation}
For later reference, we introduce the notations
\begin{equation}\label{54,1}
\boldsymbol{E}^e_{\parallel}=\boldsymbol{E}^e- (\boldsymbol{n}^0\otimes\boldsymbol{n}^0) \boldsymbol{E}^e,\qquad\boldsymbol{K}^e_{\parallel}=\boldsymbol{K}^e- (\boldsymbol{n}^0\otimes\boldsymbol{n}^0) \boldsymbol{K}^e\,.
\end{equation}
Then, from \eqref{54} we get
\begin{equation}\label{55}
\begin{array}{l}
    \boldsymbol{E}^e_{\parallel}=\tilde{E}^e_{\alpha\beta}\boldsymbol{a}_\alpha\otimes\boldsymbol{a}_\beta\,,
    \qquad
     \boldsymbol{n}^0\boldsymbol{E}^e=\boldsymbol{E}^{e,T}\boldsymbol{n}^0 =\tilde{E}^e_{3\alpha}\boldsymbol{a}_\alpha\,, \\
     \boldsymbol{K}^e_{\parallel}=\tilde{K}^e_{\alpha\beta}\boldsymbol{a}_\alpha\otimes\boldsymbol{a}_\beta\,,
    \qquad
     \boldsymbol{n}^0\boldsymbol{K}^e=\boldsymbol{K}^{e,T}\boldsymbol{n}^0 =\tilde{K}^e_{3\alpha }\boldsymbol{a}_\alpha\,.
     \end{array}
\end{equation}
If we denote the matrices by $\tilde{E}^e=\big(\tilde{E}^e_{ij}\big)_{3\times 3}\,$, $\tilde{K}^e=\big(\tilde{K}^e_{ij}\big)_{3\times 3}\,$, and also
\begin{equation}\label{55bis}
\begin{array}{l}
    \tilde{E}^e_{\parallel}=   \begin{bmatrix} \tilde{E}^e_{11} & \tilde{E}^e_{12}  \\
                                 \tilde{E}^e_{21} & \tilde{E}^e_{22}    \end{bmatrix}=
    \begin{bmatrix} \boldsymbol{d}_1  \cdot    \partial_1 \boldsymbol{y} \!-\!1 \,\, & \,\,
    \boldsymbol{d}_1  \cdot   \partial_2 \boldsymbol{y} \\
    \boldsymbol{d}_2  \cdot     \partial_1 \boldsymbol{y} \,\, & \,\, \boldsymbol{d}_2  \cdot   \partial_2 \boldsymbol{y}\!-\!1
    \end{bmatrix}, \vspace{4pt}\\
    \big(\tilde{E}^{e,T}n^0\big)= \begin{bmatrix}       \tilde{E}^e_{31} & \tilde{E}^e_{32}    \end{bmatrix}=
     \begin{bmatrix}    \boldsymbol{d}_3 \! \cdot    \partial_1 \boldsymbol{y} \,\, & \,\,
     \boldsymbol{d}_3 \! \cdot   \partial_2 \boldsymbol{y}     \end{bmatrix},
    \vspace{4pt}\\
     \tilde{K}^e_{\parallel}= \begin{bmatrix} \tilde{K}^e_{11} & \tilde{K}^e_{12}  \\
                                 \tilde{K}^e_{21} & \tilde{K}^e_{22}    \end{bmatrix} =
    \begin{bmatrix} \boldsymbol{d}_3 \! \cdot    \partial_1 \boldsymbol{d}_2 \!-\! \boldsymbol{a}_3 \! \cdot    \partial_1 \boldsymbol{a}_2 \,\, &
    \,\, \boldsymbol{d}_3 \! \cdot    \partial_2 \boldsymbol{d}_2 \!-\! \boldsymbol{a}_3 \! \cdot    \partial_2 \boldsymbol{a}_2   \\
    \boldsymbol{d}_1 \! \cdot    \partial_1 \boldsymbol{d}_3 \!-\! \boldsymbol{a}_1 \! \cdot    \partial_1 \boldsymbol{a}_3 \,\,   & \,\, \boldsymbol{d}_1 \! \cdot    \partial_2 \boldsymbol{d}_3  \!-\! \boldsymbol{a}_1 \! \cdot    \partial_2 \boldsymbol{a}_3   \end{bmatrix},
                                \vspace{4pt} \\
    \big(\tilde{K}^{e,T}n^0\big)= \begin{bmatrix}       \tilde{K}^e_{31} & \tilde{K}^e_{32}    \end{bmatrix}=
    \begin{bmatrix}   \boldsymbol{d}_2  \cdot    \partial_1 \boldsymbol{d}_1 \!-\! \boldsymbol{a}_2  \cdot    \partial_1 \boldsymbol{a}_1  \,\, & \,\, \boldsymbol{d}_2  \cdot    \partial_2 \boldsymbol{d}_1 \!-\! \boldsymbol{a}_2  \cdot    \partial_2 \boldsymbol{a}_1  \end{bmatrix},
     \end{array}
\end{equation}
 then the relations \eqref{54} and \eqref{55} can be written in matrix form
\begin{equation}\label{55,1}
\begin{array}{l}
    \tilde{E}^e=\begin{bmatrix} \tilde{E}^e_{11} & \tilde{E}^e_{12} &  0 \\
                                 \tilde{E}^e_{21} & \tilde{E}^e_{22} &  0 \\
                                 \tilde{E}^e_{31} & \tilde{E}^e_{32} & 0
    \end{bmatrix}\! = \!\! \begin{bmatrix} \big(\,\tilde{E}^e_{\parallel}\,\big)_{2\times 2} & 0_{2\times 1}  \vspace{10pt} \\
                                 \big(\tilde{E}^{e,T}\!n^0\big)_{1\times 2} & 0
    \end{bmatrix}_{3\times 3}\!\! \!\!=\!\! \begin{bmatrix} \boldsymbol{d}_1 \! \cdot    \partial_1 \boldsymbol{y} \!-\!1 & \boldsymbol{d}_1 \! \cdot   \partial_2 \boldsymbol{y} &  \!\!0 \\
                                \boldsymbol{d}_2 \! \cdot     \partial_1 \boldsymbol{y}  & \boldsymbol{d}_2 \! \cdot   \partial_2 \boldsymbol{y}\!-\!1 & \!\! 0 \\
                                 \boldsymbol{d}_3 \! \cdot    \partial_1 \boldsymbol{y}  & \boldsymbol{d}_3 \! \cdot   \partial_2 \boldsymbol{y} &  \!\!0
    \end{bmatrix}\!,\vspace{4pt}\\
    \tilde{K}^e=\begin{bmatrix} \tilde{K}^e_{11} & \tilde{K}^e_{12} &  0 \\
                                 \tilde{K}^e_{21} & \tilde{K}^e_{22} &  0 \\
                                 \tilde{K}^e_{31} & \tilde{K}^e_{32} & 0
    \end{bmatrix}  =  \begin{bmatrix} \big(\,\tilde{K}^e_{\parallel}\,\big)_{2\times 2} & 0_{2\times 1}  \vspace{10pt} \\
                                 \big(\tilde{K}^{e,T}\!n^0\big)_{1\times 2} & 0
    \end{bmatrix}_{3\times 3} \vspace{4pt}\\
    \qquad\qquad = \tilde{K}- \tilde{K}^0 = \begin{bmatrix} \boldsymbol{d}_3 \! \cdot    \partial_1 \boldsymbol{d}_2  & \boldsymbol{d}_3 \! \cdot    \partial_2 \boldsymbol{d}_2  &  0 \\
                                \boldsymbol{d}_1 \! \cdot    \partial_1 \boldsymbol{d}_3   & \boldsymbol{d}_1 \! \cdot    \partial_2 \boldsymbol{d}_3  &  0 \\
                                \boldsymbol{d}_2 \! \cdot    \partial_1 \boldsymbol{d}_1   & \boldsymbol{d}_2 \! \cdot    \partial_2 \boldsymbol{d}_1  & 0
    \end{bmatrix} -
    \begin{bmatrix} \boldsymbol{a}_3 \! \cdot    \partial_1 \boldsymbol{a}_2  & \boldsymbol{a}_3 \! \cdot    \partial_2 \boldsymbol{a}_2  &  0 \\
                                \boldsymbol{a}_1 \! \cdot    \partial_1 \boldsymbol{a}_3   & \boldsymbol{a}_1 \! \cdot    \partial_2 \boldsymbol{a}_3  &  0 \\
                                \boldsymbol{a}_2 \! \cdot    \partial_1 \boldsymbol{a}_1   & \boldsymbol{a}_2 \! \cdot    \partial_2 \boldsymbol{a}_1  & 0
    \end{bmatrix}
         \end{array}
\end{equation}
These expressions are completely similar to the strain measures for planar--shells introduced in \cite{Neff_plate04_cmt,Neff_plate07_m3as}.

Let us discuss next some important classes of elastic shells.

\subsection{Isotropic shells}

In the resultant 6-parameter theory of shells, the strain energy density for isotropic shells has been presented in various forms. The simplest expression of $W(\boldsymbol{E}^e,\boldsymbol{K}^e)$ has been proposed in the papers \cite{Pietraszkiewicz-book04,Pietraszkiewicz10} in the form
\begin{equation}\label{56}
   \begin{array}{l}
    2W(\boldsymbol{E}^e,\boldsymbol{K}^e)= \,\,\,C\big[\,\nu \,(\mathrm{tr} \boldsymbol{E}^e_{\parallel})^2 +(1-\nu)\, \mathrm{tr}(\boldsymbol{E}^{e,T}_{\parallel}  \boldsymbol{E}^e_{\parallel} )\big]  + \alpha_{s\,}C(1-\nu) \, \boldsymbol{n}^0 \boldsymbol{E}^e \boldsymbol{E}^{e,T}  \boldsymbol{n}^0 \\
    \qquad\qquad\qquad\,\, +\,D \big[\,\nu\,(\mathrm{tr} \boldsymbol{K}^e_{\parallel})^2 + (1-\nu)\, \mathrm{tr}(\boldsymbol{K}^{e,T}_{\parallel}  \boldsymbol{K}^e_{\parallel} )\big]  + \alpha_{t\,}D(1-\nu) \,    \boldsymbol{n}^0 \boldsymbol{K}^e \boldsymbol{K}^{e,T}  \boldsymbol{n}^0,
\end{array}
\end{equation}
where $C=\frac{E\,h}{1-\nu^2}\,$ is the stretching (in-plane) stiffness of the shell, $D=\frac{E\,h^3}{12(1-\nu^2)}\,$ is the bending stiffness, $h$ is the thickness of the shell, and $\alpha_s\,$, $\alpha_t$ are two shear correction factors. Also,    $E $ and $ \nu$ denote the Young modulus and  Poisson ratio of the isotropic and homogeneous material. By the numerical treatment of non-linear shell problems, the values of the shear correction factors have been set to  $\alpha_s=5/6$, $\alpha_t=7/10$ in \cite{Pietraszkiewicz10}. The value $\alpha_s=5/6$ is a classical suggestion, which has been previously deduced analytically by Reissner in the case of plates \cite{Reissner45,Naghdi72}.
Also, the value $\,\alpha_t=7/10\,$ was proposed earlier in \cite[see p.78]{Pietraszkiewicz-book79} and has been suggested in the work \cite{Pietraszkiewicz79}.
However, the discussion concerning the possible values of  shear correction factors for shells is long and controversial in the literature \cite{Naghdi72,Naghdi-Rubin95}.

With the help of the matrices \eqref{55bis}, we can express the strain energy density \eqref{56}   in the alternative form
\begin{equation}\label{57}
\begin{array}{l}
    2W(\boldsymbol{E}^e,\boldsymbol{K}^e)= \\
    =C(1\!-\!\nu)\big(\|\text{dev}_2\,\text{sym}\, \tilde{E}^e_{\parallel}\,\|^2\!+ \|\text{skew}\, \tilde{E}^e_{\parallel}\,\|^2\big)\!+
    C \dfrac{1\!+\!\nu}{2}\big(\text{tr}\, \tilde{E}^e_{\parallel}\big)^2
         + \alpha_sC(1\!-\!\nu)\,\|\tilde{E}^{e,T}n^0  \|^2 \vspace{2pt}\\
    +D(1\!-\!\nu)\big(\|\text{dev}_2\,\text{sym}\, \tilde{K}^e_{\parallel}\,\|^2\!+ \|\text{skew}\, \tilde{K}^e_{\parallel}\,\|^2\big)\!
         + D \dfrac{1\!+\!\nu}{2}\big(\text{tr}\, \tilde{K}^e_{\parallel}\big)^2 + \alpha_tD(1\!-\!\nu)\,\|\tilde{K}^{e,T}n^0   \|^2\!,
\end{array}
\end{equation}
where $\,\text{sym}\, X= \frac{1}{2}\big(X+X^T\big)$ is the symmetric part, $\,\text{skew}\, X= \frac{1}{2}\big(X-X^T\big)$ is the skew-symmetric part,
and $\,\,\text{dev}_2\, X= X-\frac{1}{2}\big(\text{tr}\,X\big)\id_2\,\,$ is the deviatoric part of any $2\times 2$ matrix  $X $. The coefficients in \eqref{57} are expressed in terms of the Lam\'e constants of the material $\lambda$ and $\mu$  by the relations
\begin{equation*}% \label{57,1}
    C \dfrac{1\!+\!\nu}{2}\,=h\,\dfrac{\mu(2\mu\!+\!3\lambda)}{2\mu+\lambda}\,,\quad\! C(1\!-\!\nu)=2\mu h,\quad\!
    D \dfrac{1\!+\!\nu}{2}\,=\dfrac{h^3}{12}\,\,\dfrac{\mu(2\mu\!+\!3\lambda)}{2\mu+\lambda}\,,\quad\!
    D(1\!-\!\nu)=\dfrac{\mu h^3}{6}\,.
\end{equation*}
Then, we obtain  that the given quadratic form \eqref{57} is positive definite if and only if the coefficients $E$ and $\nu$ satisfy the inequalities
\begin{equation}\label{57,1}
    E>0,\qquad -1<\nu<\dfrac{1}{2}\,\,.
\end{equation}
In terms of  the  Lam\'e moduli of the material,   the  inequalities \eqref{57,1} are equivalent to
$$\mu>0,\qquad 2\mu+3\lambda>0.$$
These conditions are guaranteed by the positive definiteness of the 3D quadratic elastic strain energy for isotropic materials.
Thus, we find that the strain energy $W$ is convex and satisfies the coercivity condition \eqref{26bis}, so that the hypotheses of Theorem \ref{th1} are fulfilled. Applying Theorem \ref{th1} we obtain (under suitable assumptions on the given load and boundary data, and the reference configuration $( \boldsymbol{y}^0,\boldsymbol{Q}^{0})$) the existence of minimizers for isotropic shells with strain energy density in the form \eqref{56}.\smallskip

In \cite{Eremeyev06}, Eremeyev and Pietraszkiewicz have proposed a more general form of the strain energy density, namely
\begin{equation}\label{58}
    \begin{array}{l}
   2 W(\boldsymbol{E}^e,\boldsymbol{K}^e)=  \alpha_1\big(\mathrm{tr}  \boldsymbol{E}^e_{\parallel}\big)^2 +\alpha_2  \mathrm{tr} \big(\boldsymbol{E}^e_{\parallel}\big)^2    + \alpha_3 \mathrm{tr}\big(\boldsymbol{E}^{e,T}_{\parallel}  \boldsymbol{E}^e_{\parallel} \big)  + \alpha_4     \boldsymbol{n}^0 \boldsymbol{E}^e \boldsymbol{E}^{e,T}  \boldsymbol{n}^0 \\
    \qquad\qquad\qquad + \beta_1\big(\mathrm{tr}  \boldsymbol{K}^e_{\parallel}\big)^2 +\beta_2  \mathrm{tr} \big(\boldsymbol{K}^e_{\parallel}\big)^2    + \beta_3 \mathrm{tr}\big(\boldsymbol{K}^{e,T}_{\parallel} \boldsymbol{K}^e_{\parallel} \big)  + \beta_4     \boldsymbol{n}^0 \boldsymbol{K}^e \boldsymbol{K}^{e,T}  \boldsymbol{n}^0.
\end{array}
\end{equation}
The eight coefficients $\alpha_k\,$, $\beta_k$ ($k=1,2,3,4$) can depend in general on the structure curvature tensor  $ \boldsymbol{K}^0=\text{axl}\big(\partial_\alpha \boldsymbol{Q}^{0}\boldsymbol{Q}^{0,T}\big) \otimes \boldsymbol{a}^\alpha$ of the reference configuration. For the sake of simplicity, we assume in our discussion that the coefficients $\alpha_k$ and $\beta_k$ are constant.
We can decompose the strain energy density  \eqref{58}  in the  in-plane part  $W_{\text{plane}}(\boldsymbol{E}^e)$ and the curvature part $W_{\text{curv}}(\boldsymbol{K}^e)$ and write their expressions using the matrices of components \eqref{55bis} in  form
\begin{equation}\label{59}
    \begin{array}{c}
    W(\boldsymbol{E}^e,\boldsymbol{K}^e)= W_{\text{plane}}(\boldsymbol{E}^e\big)+ W_{\text{curv}}(\boldsymbol{K}^e\big)\,,
\end{array}
\end{equation}
\begin{equation*}
    \begin{array}{c}
    2W_{\text{plane}}(\boldsymbol{E}^e)= (\alpha_2\!+\!\alpha_3) \| \text{sym}\, \tilde{E}^e_{\parallel}\|^2\!+ (\alpha_3\!-\!\alpha_2)\|\text{skew}\, \tilde{E}^e_{\parallel}\|^2\!+ \alpha_1\big(\mathrm{tr}  \tilde{E}^e_{\parallel}\big)^2 +  \alpha_4\|\tilde{E}^{e,T}n^0  \|^2, \vspace{4pt}\\
    2 W_{\text{curv}}(\boldsymbol{K}^e)= (\beta_2\!+\!\beta_3) \| \text{sym}\, \tilde{K}^e_{\parallel}\|^2\!+ (\beta_3\!-\!\beta_2)\|\text{skew}\, \tilde{K}^e_{\parallel}\|^2\!+ \beta_1\big(\mathrm{tr}  \tilde{K}^e_{\parallel}\big)^2
    + \beta_4\|\tilde{K}^{e,T}n^0  \|^2.
    \end{array}
\end{equation*}
The in-plane part of the energy density \eqref{59} can be written equivalently as
\begin{equation}\label{59,1}
 \begin{array}{l}
    2W_{\text{plane}}(\boldsymbol{E}^e)= \underbrace{(\alpha_2+\alpha_3)\, \|\,\text{dev}_2 \,\text{sym}\, \tilde{E}^e_{\parallel}\,\|^2 }_{\text{in-plane shear--stretch energy}}  \,+\, \underbrace{(\alpha_3-\alpha_2)\,\|\,\text{skew}\, \tilde{E}^e_{\parallel}\,\|^2\,}_{\text{in--plane drill rotation energy}}\vspace{4pt}\\
    \qquad\qquad\quad\,\,\, +\,\underbrace{\Big(\alpha_1+\dfrac{\alpha_2 + \alpha_3}{2}\Big)\big(\mathrm{tr} \, \tilde{E}^e_{\parallel}\big)^2}_{\text{in-plane elongational stretch energy}}  \,+\,\underbrace{\, \alpha_4\,\|\,\tilde{E}^{e,T}n^0 \, \|^2}_{\text{transverse shear energy}}.
    \end{array}
\end{equation}
The above forms of the strain energy $W$ are expressed in terms of the components of the tensors $\boldsymbol{E}^e$ and $\boldsymbol{K}^e$ in the basis $\{\boldsymbol{a}_i\otimes \boldsymbol{a}_\alpha\}\,$, i.e. in terms of the elements of the matrices \eqref{55bis}.
Denoting with $\,\,\mu_c^{\text{drill}}\,$ the coefficient $\,(\alpha_3-\alpha_2)\,$  in   \eqref{59,1}, we remark that the   term
\begin{equation}\label{59,2}
    \mu_c^{\text{drill}}\,\|\,\text{skew}\, \tilde{E}^e_{\parallel}\,\|^2\,,\qquad \text{with}\qquad \mu_c^{\text{drill}}\,:=\,\alpha_3-\alpha_2\,,
\end{equation}
describes the quadratic in-plane drill rotation energy of the shell. We call the coefficient $\,\,\mu_c^{\text{drill}}\,$  the \emph{linear in-plane rotational couple modulus}, in analogy to the Cosserat couple modulus in the three-dimensional Cosserat theory \cite{Neff_zamm06}.
\begin{remark}
The planar isotropic Cosserat shells have been investigated also in \cite{Neff_plate04_cmt,Neff_plate07_m3as}, using a model derived directly from the 3D equations of Cosserat elasticity. We mention that the expressions \eqref{59}, \eqref{59,1} of the strain energy density are essentially the same as the strain energy of the Cosserat model for planar shells \cite{Neff_plate04_cmt}. By comparing these two approaches (6-parameter resultant shells and Cosserat model) we deduce the following identification of the constitutive coefficients $\alpha_1\,,...,\alpha_4$
\begin{equation}\label{59,3}
    \alpha_1=h\,\dfrac{2\mu\lambda}{2\mu+\lambda}\,,\quad \alpha_2=h(\mu-\mu_c),\quad \alpha_3=h(\mu+\mu_c),\quad \alpha_4=\kappa\, h (\mu+\mu_c),
\end{equation}
where $\,\mu_c\,$ is the Cosserat couple modulus of the 3D continuum, and $\,\kappa\,$ is a formal shear correction factor. From \eqref{59,2}, \eqref{59,3} we observe that
\begin{equation}\label{59,4}
    \mu_c^{\mathrm{drill}}\,=\,\alpha_3-\alpha_2\,=\,2h\,\mu_c\,,
\end{equation}
which means that the in-plane rotational couple modulus $\,\mu_c^{\mathrm{drill}}\,$ of the Cosserat shell model is determined by the Cosserat couple modulus $\,\mu_c\,$ of the 3D Cosserat material.

The relations \eqref{59,3} are similar to the corresponding relations in the linear theory of micropolar plates, see \cite[Eqs.(45)]{Altenbach-Erem09}. From a mathematical viewpoint, the difference between the two sets of relations consists in the notations used and the value of the shear correction factor.
\end{remark}

Looking at \eqref{59} and \eqref{59,1} we observe that  the quadratic form $ W(\boldsymbol{E}^e,\boldsymbol{K}^e) $ is positive definite  if and only if the coefficients verify the conditions
\begin{equation}\label{60}
    \begin{array}{l}
    2\alpha_1+\alpha_2+\alpha_3>0,\quad \alpha_2+\alpha_3>0,\quad \alpha_3-\alpha_2>0, \quad\alpha_4>0,\\
    2\beta_1+\beta_2+\beta_3>0,\quad \beta_2+\beta_3>0,\quad \beta_3-\beta_2>0,\quad  \beta_4>0,
\end{array}
\end{equation}
Provided that the conditions \eqref{60} are satisfied, the strain energy function $W(\boldsymbol{E}^e,\boldsymbol{K}^e)$ is convex and coercive in the sense of \eqref{26bis}. By virtue of Theorem \ref{th1}, in this case the minimization problem associated to the deformation of isotropic elastic shells admits at least one solution.

\begin{remark}
The same conditions \eqref{60} have been imposed in \cite{EremeyevLebedev11} to establish existence results in the \emph{linearized} theory of micropolar (6-parameter) shells.
\end{remark}
\begin{remark}
The case $\,\mu_c^{\mathrm{drill}}=0\,$ (i.e., $\alpha_3-\alpha_2=0$) is not uniformly positive definite. However, with a slight change of the resultant shell model, one can prove the existence of minimizers using similar methods as in \cite{Neff_plate07_m3as}. A linearization of such a model  leads exactly to the Reissner kinematics with 5 degrees of freedom \cite{Neff_plate07_m3as}, where the in-plane drill rotation is absent. The physical meaning of the in-plane rotational stiffness $\,\mu_c^{\mathrm{drill}}=\alpha_3-\alpha_2\,$ in the resultant shell model is not entirely clear to us.

Since only two independent rotations are required to orient a unit director field, a distinctive feature of classical plate and shell theories is a rotation field defined in terms of only \emph{two} independent degrees of freedom. Rotations about the director itself -- the so-called drill rotation -- are irrelevant and for that matter undefined in classical shell theory.
\end{remark}

\subsection{Orthotropic shells}

The constitutive equations for orthotropic shells have been presented in \cite{Eremeyev06} within the 6-parameter resultant shell  theory. The expression of the strain energy density in terms of the tensor components defined in \eqref{54} is given by
\begin{equation}\label{61}
   2W(\boldsymbol{E}^e,\boldsymbol{K}^e)= C_{\alpha\beta\gamma\delta}^E\, \tilde{E}^e_{\alpha\beta} \tilde{E}^e_{\gamma\delta} +  D_{\alpha\beta}^E\, \tilde{E}^e_{3\alpha } \tilde{E}^e_{3\beta} +
      C_{\alpha\beta\gamma\delta}^K\,   \tilde{K}^e_{\alpha\beta} \tilde{K}^e_{\gamma\delta}+      D_{\alpha\beta}^K\, \tilde{K}^e_{3\alpha } \tilde{K}^e_{3\beta}
\end{equation}
where $C_{\alpha\beta\gamma\delta}^E\,$, $C_{\alpha\beta\gamma\delta}^K\,$, $D_{\alpha\beta}^E\,$ and $D_{\alpha\beta}^K\,$ are material constants which satisfy the following symmetry relations
$$ C_{\alpha\beta\gamma\delta}^E=  C_{\gamma\delta\alpha\beta}^E\,,\quad D_{\alpha\beta}^E= D_{\beta\alpha}^E\,,\quad
      C_{\alpha\beta\gamma\delta}^K=  C_{\gamma\delta\alpha\beta}^K\,,\quad D_{\alpha\beta}^K= D_{\beta\alpha}^K\,.$$
We observe that the quadratic function \eqref{61} is coercive if and only if the following symmetric matrices are positive definite
\begin{equation}\label{62}
    \begin{bmatrix} C_{1111}^E & C_{1122}^E  &  C_{1112}^E  &   C_{1121}^E   \\
                    C_{1122}^E & C_{2222}^E  &   C_{2212}^E  &   C_{2221}^E   \\
                    C_{1112}^E  & C_{2212}^E   &  C_{1212}^E &   C_{1221}^E  \\
                    C_{1121}^E  &  C_{2221}^E   &  C_{1221}^E   &  C_{2121}^E   \end{bmatrix},
    \begin{bmatrix} C_{1111}^K & C_{1122}^K  &  C_{1112}^K  &   C_{1121}^K   \\
                    C_{1122}^K & C_{2222}^K  &   C_{2212}^K  &   C_{2221}^K   \\
                    C_{1112}^K  & C_{2212}^K   &  C_{1212}^K &   C_{1221}^K  \\
                    C_{1121}^K  &  C_{2221}^K   &  C_{1221}^K   &  C_{2121}^K   \end{bmatrix},
    \begin{bmatrix} D_{11 }^E & D_{12}^E   \\
                    D_{12}^E & D_{22}^E    \end{bmatrix},
    \begin{bmatrix} D_{11 }^K & D_{12}^K   \\
                    D_{12}^K & D_{22}^K    \end{bmatrix}\!.
\end{equation}
In the situation when the matrices \eqref{62} are positive definite, then the strain energy $W$ given by \eqref{61} satisfies the hypotheses of Theorem \ref{th1}. Then, we can use our theoretical results to derive the existence of minimizers for orthotropic shells.

\subsection{Composite layered shells}

Let us analyze the case of composite shells made of a finite number of individually homogeneous layers. According to \cite{Chroscielewski11}, the strain energy density of such type of shells can be written by means of the tensor components \eqref{54} in the form
\begin{equation}\label{63}
\begin{array}{l}
    2W(\boldsymbol{E}^e,\boldsymbol{K}^e)= A_{\alpha\beta\gamma\delta} \, \tilde{E}^e_{\alpha\beta} \tilde{E}^e_{\gamma\delta} +  D_{\alpha\beta\gamma\delta} \, \tilde{K}^e_{\alpha \beta} \tilde{K}^e_{\gamma\delta} +
      B_{\alpha\beta\gamma\delta} ( \tilde{E}^e_{\alpha\beta} \tilde{K}^e_{\gamma\delta} + \tilde{K}^e_{\alpha\beta} \tilde{E}^e_{\gamma\delta})\vspace{4pt}\\
      \qquad\qquad\qquad\quad
      +      S_{\alpha\beta} \, \tilde{E}^e_{3\alpha } \tilde{E}^e_{3\beta}
        +      G_{\alpha\beta} \,        \tilde{K}^e_{3\alpha } \tilde{K}^e_{3\beta}\,,
        \end{array}
\end{equation}
where $A_{\alpha\beta\gamma\delta}\,$, $B_{\alpha\beta\gamma\delta}\,$, $D_{\alpha\beta\gamma\delta}\,$, $S_{\alpha\beta}\,$ and $G_{\alpha\beta}\,$ are the constitutive coefficients of composite elastic shells, which have been determined in \cite{Chroscielewski11} in terms of the material/geometrical parameters of the layers. They satisfy the symmetry conditions
$$A_{\alpha\beta\gamma\delta}=  A_{\gamma\delta\alpha\beta}\,,\quad D_{\alpha\beta\gamma\delta} = D_{\gamma\delta\alpha\beta}\,,\quad
      S_{\alpha\beta}=  S_{\beta\alpha}\,,\quad G_{\alpha\beta} = G_{\beta\alpha}\,.$$
In the constitutive relation \eqref{63} one can observe a multiplicative coupling of the strain tensor $\boldsymbol{E}^e$ with the curvature tensor $\boldsymbol{K}^e$ for composite shells. Let us denote by $A$, $D$ and $B$ the $4\times 4$ matrices of material constants
\begin{equation*}
\begin{array}{c}
    A=\begin{bmatrix} A_{1111}  & A_{1122}   &  A_{1112}   &   A_{1121}    \\
                    A_{1122}  & A_{2222}   &   A_{2212}   &   A_{2221}    \\
                    A_{1112}   & A_{2212}    &  A_{1212}  &   A_{1221}   \\
                    A_{1121}   &  A_{2221}    &  A_{1221}    &  A_{2121}    \end{bmatrix},\qquad
    D=\begin{bmatrix} D_{1111}  & D_{1122}   &  D_{1112}   &   D_{1121}    \\
                    D_{1122}  & D_{2222}   &   D_{2212}   &   D_{2221}    \\
                    D_{1112}   & D_{2212}    &  D_{1212}  &   D_{1221}   \\
                    D_{1121}   &  D_{2221}    &  D_{1221}    &  D_{2121}    \end{bmatrix},\\
    B=\begin{bmatrix} B_{1111}  & B_{1122}   &  B_{1112}   &   B_{1121}    \\
                    B_{2211}  & B_{2222}   &   B_{2212}   &   B_{2221}    \\
                    B_{1211}   & B_{1222}    &  B_{1212}  &   B_{1221}   \\
                    B_{2111}   &  B_{2122}    &  B_{2112}    &  B_{2121}    \end{bmatrix}.
\end{array}
\end{equation*}
One can show that the necessary and sufficient condition for the coercivity of the strain energy function \eqref{63} is that the following matrices are positive definite
\begin{equation*}
    \mathbb{C}=\begin{bmatrix} A_{4\times 4 }  & B_{4\times 4 }   \\
                    B_{4\times 4 } & D_{4\times 4 }    \end{bmatrix}_{8\times 8}\,,\quad
                    \mathbb{S}= \begin{bmatrix} S_{11 }  & S_{12}    \\
                    S_{12}  & S_{22}     \end{bmatrix}_{2\times 2}\,,\quad
                    \mathbb{G}= \begin{bmatrix} G_{11 }  & G_{12}    \\
                    G_{12}  & G_{22}     \end{bmatrix}_{2\times 2}\,.
\end{equation*}
With these notations, one may write the strain energy density \eqref{63} in the matrix form
\begin{equation*}
\begin{array}{l}
   2W(\boldsymbol{E}^e,\boldsymbol{K}^e)= V\,\mathbb{C}\,V^T+ \big(\,\tilde{E}^e_{31}\,,\,\tilde{E}^e_{32}\big)\, \mathbb{S}\, \big(\,\tilde{E}^e_{31}\,,\,\tilde{E}^e_{32}\big)^T + \big(\,\tilde{K}^e_{31}\,,\,\tilde{K}^e_{32}\big)\, \mathbb{G}\, \big(\,\tilde{K}^e_{31}\,,\,\tilde{K}^e_{32}\big)^T, \vspace{4pt}\\
   \qquad\quad\text{with}\qquad V= \big(\,\tilde{E}^e_{11}\,,\,\tilde{E}^e_{22}\,,\, \tilde{E}^e_{12}\,,\,\tilde{E}^e_{21}\,,\,\tilde{K}^e_{11}\,,\,\tilde{K}^e_{22}\,,\, \tilde{K}^e_{12}\,,\,\tilde{K}^e_{21}\big)_{1\times 8}\,\,.
\end{array}
\end{equation*}
In conclusion, if the matrices $\mathbb{C}$, $\mathbb{S}$ and $\mathbb{G}$ are positive definite, then we can apply  Theorem \ref{th1} for the strain energy density given by \eqref{63} and prove the existence of minimizers for composite layered shells.

\begin{remark}
The results and conclusions presented above are obviously valid also in the case of \emph{plates}, i.e. when the reference base surface $S^0$ is \emph{planar}. However, many of the formulas  for general shells can be significantly simplified in the case of plates, since   the 3 orthonormal bases
$\{\boldsymbol{a}_1,\boldsymbol{a}_2,\boldsymbol{n}^0\}\,$, $\{\boldsymbol{d}^0_1,\boldsymbol{d}^0_2,\boldsymbol{d}^0_3\}$ and $\{\boldsymbol{e}_1,\boldsymbol{e}_2,\boldsymbol{e}_3\}$ can be considered identical.

The corresponding existence results for 6-parameter geometrically non-linear plates (planar shells) has been presented in \cite{Birsan-Neff-Plates} for isotropic or anisotropic materials, and in \cite{Birsan-Neff-AnnRom12} for composite planar--shells. We mention that, in the case of isotropic plates, the existence theorem can be obtained from the more general results concerning Cosserat planar--shells presented by the second author in \cite{Neff_plate04_cmt,Neff_plate07_m3as}.
\end{remark}

In a forthcoming contribution we will extend our existence results to the 6-parameter resultant shell model with physically non-linear behavior and show the invertibility of the reconstructed deformation gradient $\bar{F}$.

\bigskip\bigskip
\small{\textbf{Acknowledgements.}
The first author (M.B.) is supported by the german state  grant: ``Programm des Bundes und der L\"ander f\"ur bessere Studienbedingungen und mehr Qualit\"at in der Lehre''. We thank our many friends who have made substantial comments on a preliminary version of the paper.

%************************************************
%************************************************
%***************************************************

\bibliographystyle{plain} %plain
{\footnotesize
\bibliography{literatur_Birsan}
%\bibliography{/home/patrizio/pasadena/plasticity/literatur1}
%\bibliography{/home/neff/pasadena/plasticity/literatur1}
%\bibliography{/Users/patrizio/Library/texmf/literatur1}
%        \bibliography{/home/nesenenko/publications/GradientPlasticity/literatur1,
%\home/nesenenko/publications/GradientPlasticity/literaturSN}
%\bibliography{/home/nesenenko/publications/GradientPlasticity/literaturSN}
}

%***************************************************

\end{document}